\documentclass[11pt]{article}
\usepackage{amsmath}
\usepackage{amscd}
\usepackage{epsfig}
\usepackage[latin1]{inputenc}

\input{amssym}

\newcommand{\color}[6]{}
\newcommand\QQ{\hbox{I\kern-.53em\hbox{Q}}}
\newcommand\qed{\hfill$\sqcap\kern-8.0pt\hbox{$\sqcup$}$}
\newcommand\NN{\hbox{I\kern-.2em\hbox{N}}}
\newcommand\RR{\hbox{I\kern-.2em\hbox{R}}}
\newcommand\sRR{{\sl \hbox{I\kern-.2em\hbox{R}}}}
\newcommand{\PP}{{\bf P}^k}
\newcommand{\pp}{{\bf P}^1}
\newcommand{\fL}{\{f_{\la}\}_{\la \in X}}
\newcommand\ZZ{{{\rm Z}\kern-.28em{\rm Z}}}
\newcommand\proof{\noindent{\em{Proof}.\ }}
\newcommand{\Hd}{{\cal H}_d}
\newtheorem{theo}{Theorem}[section]
\newtheorem{prop}[theo]{Proposition}
\newtheorem{lem}[theo]{Lemma}
\newtheorem{rem}[theo]{Remark} 
\newtheorem{cor}[theo]{Corollary}
\newtheorem{defi}[theo]{Definition}

\newcommand{\la}{\lambda}

\numberwithin{equation}{section}

\begin{document}

\date{}
\title {Bifurcation currents in holomorphic dynamics on $\PP$}

\vskip0.5cm

\author{Giovanni Bassanelli and Fran\c{c}ois Berteloot\thanks{The authors thank the Laboratoire E. Picard de Toulouse 
and the Dipartimento di Matematica di Parma for the reciprocal
 kind hospitality
during the preparation of the paper. The first author was partially supported by GNSAGA of INdAM, the second one
by the Universit\' e Italo-Fran\c caise.}}

\vskip0.5cm

\maketitle

\section*{Introduction}

Potential theory has been introduced in one-dimensional rational dynamics by Brolin and Tortrat
(\cite{Bro}, \cite{Tor}) but   it   does not play  a central role there. In higher dimension however, as the classical tools are not any longer efficient, 
pluri-potential theory has revealed itself to be essential. The fundamental works of Hubbard-Papadopol, Fornaess-Sibony,  
Briend-Duval,   Bedford-Smillie   (see \cite{Sb} for precise references) enlighten the remarkable  effectiveness   of 
pluri-potential theory in holomorphic dynamics 
on $\PP$ or ${\bf C}^{k}$. It is therefore tempting to 
study the parameter spaces in a similar way. More precisely, one would like to relate the bifurcations of a holomorphic family $\fL$ of 
endomorphisms of $\PP$ to the powers of a certain current on the parameter space $X$. \\
Let us recall that in dimension $k=1$, a bifurcation is said to occur at some point $\lambda_0\in X$ if the Julia set of $f_{\la}$ does not 
move continuously around $\la_0$. The famous work of Ma\~n\'e-Sad-Sullivan \cite{MSS}, which is based on the $\la$-lemma and the 
Fatou-Cremer-Sullivan classification, relates the bifurcations with the instability of the critical orbits. It 
also   asserts   that the 
bifurcations concentrate on the complement of an open dense subset of $X$ (for the quadratic family $\{z^2+\la\}_{ \la \in X={\bf C}}$
the bifurcation locus is precisely the boundary of the Mandelbrot set).\\

A seminal idea towards the application of  potential theory to the study of bifurcations is due to 
 Przytycki,   who raised the following problem in the final remarks of his paper.
 {\em Problem   \cite{Prz}: understand the connections between  Lyapunov characteristic exponents and potential theory
 for rational mappings.}\\
To support his point of view, Przytycki also analysed the following formula for a polynomial $p$ of degree $d$ on ${\bf C}$ 
(see also \cite{Mann}):
\begin{equation}\label{i1} L(p)=\sum_j G_p(c_j) + \log d \end{equation}
where $L(p)$ is the  Lyapunov exponent of $p$  with respect to its equilibrium measure,
$G_p$ its Green function  and $c_j$ are the critical points of $p$.
More recently, DeMarco (\cite{Demar}, \cite{Demar2})  has obtained a generalization of this formula to rational maps of $\pp$.
 She also used her formula to show that, for a holomorphic family $\fL$, 
the current $dd^cL(f_{\la})$ is supported by the bifurcation locus.\\

The results of the present paper deal with Przytycki problem. Our first goal is to show that $dd^cL(f_\la)$ is a reasonable bifurcation current
in any dimension. To this purpose, we prove the following theorem
in Section \ref{SMvH}.
\medskip

\noindent {\bf Theorem \ref{a}} {\it Given a holomorphic family $\fL$ of 
endomorphisms of $\PP$, the function $L(f_\la)$, defined as the sum of  Lyapunov exponents of $f_\la$ for its Green measure, 
is pluriharmonic 
on $X$ if the repulsive cycles of $f_\la$ move holomorphically on $X$.}
 \medskip

Let us mention here that 
all what we need to know about $L$ in the paper is that $L(F)=\int_{{\bf C}^{k+1}} \log \vert det F'\vert \mu_F$ where $F$ is a lift of 
$f$ and $\mu_F$   is   its Green measure. \\
With the goal of analysing the support of $dd^cL(f_{\la})$ and its powers, we then generalize formula \eqref{i1} to
 endomorphisms of $\PP$. This is done in a very natural way by using an 
integration by part
 formula on a suitable line bundle. We obtain the following 

\medskip

\noindent {\bf Formula} (see Theorem \ref{o}). 
$L(f)=\sum_{j=0}^{k-1}  \int_{\PP}g_F \ (dd^cg_F+\omega)^j \wedge \omega^{k-j-1}\wedge [C_f]   -$
$$-\log d+\int_{{\bf P}^k}\log||J_F||_0\omega^k-(k+1)(d-1)
\sum_{j=0}^k\int_{{\bf P}^k}g_F \  (dd^cg_F+\omega)^j \wedge \omega^{k-j}.$$

For a holomorphic family $\fL$, the above formula allows us to compute the bifurcation current 
$dd^c L(f_\la)$. We get the following synthetic statement:

\medskip

\noindent {\bf Theorem} (see Corollary \ref{v})
\begin{equation}\label{i2}dd^cL(f_{\la})=p_*((dd^cg_{F_{\la}}+\omega)^k\wedge[C_X]),\end{equation}
{\it and on} $X \times {\bf P}^k$
\begin{equation}\label{i3}(dd^cg_{F_{\la}}+\omega)^{k+1}=0.\end{equation}

\noindent In these formulas   the operator   $dd^c$ is acting on $X\times \PP$, $p$ is the canonical projection from 
$X\times \PP$ to $X$, $g_{F_\la}$ is the Green function of $f_\la$ on $\PP$   associated   to the lift $F_{\la}$ and 
  \label{NB 1} $C_X$ is the hypersurface of $X \times \PP$ defined by the equation $det F'_{\la}(z)=0$  .

It is worth emphasize that there is a certain interaction between formulas   \eqref{i2} and \eqref{i3},   this may be seen 
in the example in Subsection \ref{S7b} and in the Appendix (see formula \eqref{ad}). 
Moreover, formula \eqref{i3} is formally equivalent to the equation of geodesics on
 the space of K\" ahler metrics on $\PP$;   this leads to some examples of such geodesics,   as discussed in 
Subsection \ref{S7a}.
These results are established in Section \ref{SF2} but   we treat the case of dimension $k=1$ separately in Section \ref{SF1},
since it   is technically less involved and may help the reader to a better understanding.
Let us also stress that in the one-dimensional case we   get   several explicit formulas for $L(f)$ (see Theorem \ref{i}).
Moreover, our approach offers a much simpler proof of DeMarco's formula. The equivalence between DeMarco's formula and 
ours is a consequence of the following identity which may be of independent interest (see Proposition \ref{z}):
$$\int_{{\bf P}^1}g_F(\mu_f+\omega)=\frac{1}{2}(\log|Res(F)|-1).$$

 Sections 5 and 6 are devoted to the study of bifurcations for holomorphic families
$\fL$ of rational maps of $\pp$. For a holomorphic family $\fL$ of endomorphisms of
$\pp$ we obtain a geometrical description of the support of the bifurcation currents $(dd^cL(f_{\la}))^p$ by means of
certain complex hypersurfaces. 

For $\theta \in {\bf R} \setminus {\bf Z}$, the  set of all $\la \in X$
such that $f_{\la}$ 
 has a periodic point of period $n$ and multiplier $e^{2i\pi\theta}$ is generically a complex hypersurface of $X$, denoted by 
$Per(X, n, e^{2i\pi \theta})$. Therefore,   for $n$ fixed  ,   $\bigcup_{\theta}Per(X, n, e^{2i\pi \theta})$ can be thought as a real
hypersurface foliated
by complex hypersurfaces. The union ${\cal Z}_1(X)$ of all these hypersurfaces is dense in the support of the
 the bifurcation current $dd^cL(f_{\la})$:

\medskip 

\noindent {\bf Theorem}\ \   $\overline{{\cal Z}_1(X)}=$ Supp$\big(dd^cL(f_{\la})\big)$.
\medskip

\noindent We included this geometrical characterization of the bifurcation locus in the statement of {\bf Theorem \ref{zb}}
which may be interpreted as treating a substantial part of  Ma\~n\'e-Sad-Sullivan theory by 
potential-theoretic methods. The proof exploits the links between the vanishing of $dd^cL$, the motion of repulsive cycles
 and the stability of critical orbits. We point out that these links are revealed by formula \eqref{i2} and  the above quoted  Theorem \ref{a}.

For the powers of $dd^cL(f_{\la})$ the geometry is more involved. 
Taking all possible intersections between $p $ of the above complex hypersurfaces one gets a large family of 
codimension $p$ subvarieties
of $X$; the union of this family, denoted by ${\cal Z}_p(X)$, satisfies:

\medskip

\noindent{\bf Theorem \ref{zf}}\ \  {\it For any} 
$1 \le p \le dim_{\bf C}X,\ Supp\big((dd^cL(f_{\la})^p\big) \subset \overline {{\cal Z}_p(X)}.$

\medskip

\noindent These results have some significant consequences
 considering  the family $\Hd(\pp)$ of all the rational maps of degree $d$. First of all one may show 
that:

\medskip

\noindent {\bf Proposition \ref{aq}} {\it For $1 \le p \le 2d+1$ the bifurcation current $(dd^cL)^p$ has finite mass on} $\Hd(\pp)$.

\medskip

\noindent This implies 
 that $(dd^cL)^{2d-2}$ induces a measure $\mu$ of finite mass on the moduli space ${\cal M}_d$ of rational maps of degree 
$d$.
  We call it the {\it bifurcation measure} and show that its support contains all isolated Latt\`es maps.
 Using our description of $Supp\big((dd^cL)^p \big)$ we also obtain the   following fact:  

\medskip

\noindent {\bf Proposition \ref{at}} {\it Any open set of ${\cal M}_d$ intersecting the support of the bifurcation measure 
contains an uncountable set of 
chaotic mappings. }

\medskip

\noindent As a by-product of our investigation one sees that any non flexible Latt\`es map is generating quite complicated bifurcations.

\section{Preliminaries}

In all the paper $\omega$ denotes the Fubini-Study form in $\PP$ and   let   $||\ ||$ be the Hermitian norm in ${\bf C}^{k+1}.$
\subsection{The spaces $\Hd({\bf C}^{k+1})$ and $\Hd(\PP)$}\label{sssp}
Every holomorphic endomorphism $f$ of $\PP$ has a lift $F:{\bf C}^{k+1} \to {\bf C}^{k+1}$  ,   that is:
 a homogeneous, non-degenerate, polynomial map such that $\pi \circ F=f \circ \pi$, where
\linebreak ${\bf C}^{k+1}\setminus \{0\} \stackrel{\pi}{\rightarrow} \PP$ is the canonical projection. 
The degree $d$ of $F$ is, by definition, the
{\em algebraic degree} of $f$, while $d^k$ is the topological degree of $f$. 

In the following it will always be assumed that
$d \geq 2$.

Since a homogeneous polynomial of degree $d$ in $k+1$
 variables depends on
$\frac{(d+k)!}{d!k!}$ coefficients, such a lift $F=(F_0, \dots, F_k)$ can be identified with an element of ${\bf C}^{N+1}$,
where $N=(k+1)\frac{(d+k)!}{d!k!}-1$. With this identification, the space  of all homogeneous, non-degenerate,
polynomial maps of degree $d$ on ${\bf C}^{k+1}$  is an open subset of ${\bf C}^{N+1}$
 which we denote by $\Hd ({\bf C}^{k+1})$.
 Denoting again by ${\bf C}^{N+1}\setminus \{0\} \stackrel{\pi}{\rightarrow} {\bf P}^N$ the canonical
projection, we get $\pi(F)=f$ and   that   $\pi(\Hd  ({\bf C}^{k+1}))$ is the space
 of all holomorphic endomorphisms of $\PP$ of degree $d$  ,   which we denote by $\Hd(\PP)$.

The following Proposition shows that the complement  of $\Hd ({\bf C}^{k+1})$ in ${\bf C}^{N+1}$ is an
  irreducible   complex hypersurface $\tilde \Sigma_d= \{Res =0\}$.  Thus the projective hypersurface
$\Sigma_d:=\pi(\tilde \Sigma_d)$ is the complement of $\Hd(\PP)$ in ${\bf P}^N$.

\begin{prop}\label{x}  Let  $F_0, \dots, F_k,$  be homogeneous polynomials of 
degree $d$, 
in $k+1$ complex variables. 
There exists a unique polynomial $Res(F_0, \dots, F_k)$  in the coefficients 
of $F_0, \dots, F_k,$ 
which is homogeneous of degree
 $(k+1)d^k$, irreducible,   and
such that
\begin{itemize}
\item[(i)]  $Res(F_0, \dots, F_k)=0$ if and only if $F=(F_0, \dots, F_k):
{\bf C}^{k+1} \to {\bf 
C}^{k+1}$ is degenerate, 

\item[(ii)] $Res(z_0^d, \dots, z_k^d)=1.$
\end{itemize}
\end{prop}

\proof See \cite{GKZ} p. 427 and p. 105.

\subsection{Green functions}

To any $F\in \Hd({\bf C}^{k+1})$ it is associated a {\it Green function} $G_F$ 
defined by 
$$G_F:=\lim_n d^{-n}\log\Vert F^n(z)\Vert.$$
Let us stress that $G_F$ is the limit of a sequence   $\{G_{F,n}\}$   of $p.s.h.$ and 
continuous functions on $\Hd({\bf C}^{k+1})\times ({\bf C}^{k+1}\setminus\{0\})$.
The following Proposition summarizes the regularity properties of $G_F(z)$. 
The only novelty here is the H\" older-continuity in $F$.

 \begin{prop}\label{1a}
i) For any compact subset ${\cal K}$
of $\Hd({\bf C}^{k+1})$, the sequence $G_{F,n}(z)$ converges uniformly to $G_F(z)$
on ${\cal K}\times ({\bf C}^{k+1}\setminus\{0\})$. In particular $G_F(z)$ is $p.s.h$ and continuous
on $\Hd({\bf C}^{k+1}) \times 
({\bf C}^{k+1}\setminus\{0\})$. It satisfies the following homogeneity property:
$$G_F(tz)=\log\vert t\vert +G_F(z);\;\;\;\forall t\in {\bf C}^*,
\forall z\in {\bf C}^{k+1}$$
and   the   functional equation
\begin{equation}\label{3}G_F \circ F=d G_F.\end{equation}
In the definition of $G_{F}$, the norm $\Vert\;\Vert$ may be replaced by any continuous gauge function.

ii) The function $G_F(z)$ is H\" older-continuous on every compact subset 
of $\Hd({\bf C}^{k+1}) \times 
 ({\bf C}^{k+1}\setminus\{0\})$. 

\end{prop}  

\proof
i) Let $N$ be any continuous gauge function on ${\bf C}^{k+1}$.
Let $\cal K$ be a compact subset of $\Hd({\bf C}^{k+1})$ and $C>1$ be a 
constant such that
$$\frac{1}{C} \le N\left(F(z)\right)\le C;\;\;\forall F\in{\cal K},\forall 
z\in\{N=1\}.$$
Then, by homogeneity we have:
\begin{equation}\label{1}
\frac{1}{C^{1+...+d^{n-1}}}N(z)^{d^n} \le N\left(F^n(z)\right) \le
C^{1+...+d^{n-1}}N(z)^{d^n};\ \forall F\in {\cal K},
\forall z \in {\bf C}^{k+1}\setminus\{0\},\forall n \in {\bf N}. 
\end{equation}

After replacing $z$ by $F^m(z)$, taking $\log$ and dividing by $d^{n+m}$,
\eqref{1} gives:
$$\vert G_{F,n+m}(z) - G_{F,n}(z)\vert \le \frac{\log C}{d^m(d-1)};\;\;
\forall F\in {\cal K}, \forall z\in {\bf C}^{k+1}\setminus \{0\}.$$
ii) The H\" older continuity in $z$ has been established by Briend-Duval (\cite
{BD}). Inspecting the proof and taking into account the continuity of $G_F(z)$ 
in $(F,z)$, it is not hard to see that the constants might be chosen uniformly 
in $F$ (see \cite{Sb} Th\'eor\`eme 1.7.1 and Remarque 1.7.2). More precisely, 
for any compact ${\cal K}\times K \subset \Hd({\bf C}^{k+1}) \times 
({\bf C}^{k+1}\setminus\{0\})$ there are constants $C>0$ and $0<\alpha<1$ such 
that:
$$\vert G_F(z)-G_F(z') \vert \le C\Vert z-z'\Vert^{\alpha};\;\;\forall F\in 
{\cal K}, \forall z,z' \in K.$$

We will show how this property may be   transferred  to   $F$.
We may assume that $\{
G_F=0\}\subset K$ for every $F\in {\cal K}$. Let us pick $F,F_0\in {\cal K}$
and consider the gauge function $N_0:=e^{G_{F_0}}$. By the H\" older-continuity of 
$G_{F_0}$ we have
$$\vert G_{F_0}\left(F(z)\right) - G_{F_0}\left(F_0(z)\right)\vert
\le C_1\Vert F(z)-F_0(z)\Vert^{\alpha}
\le C_2\Vert F-F_0\Vert^{\alpha};\; \forall z\in K.$$
When $z \in \{N_0=1\}$ this inequality becomes
$\frac{1}{C_0} \le N_0\left(F(z)\right)\le C_0$
  where $C_0:=e^{C_2\Vert F-F_0\Vert^{\alpha}}.$
Just like for \eqref{1}, this implies
\begin{equation*}
\frac{1}{C_0^{1+...+d^{n-1}}}N_0(z)^{d^n} \le N_0\left(F^n(z)\right) \le
C_0^{1+...+d^{n-1}}N_0(z)^{d^n};\forall F\in {\cal K},
\forall z \in {\bf C}^{k+1}\setminus\{0\},\forall n \in {\bf N}. 
\end{equation*}
Taking $\log$, dividing by $d^n$ and making $n\to\infty$, this yields 
(as $G_F=\lim_n d^{-n}N_0\left(F^n(z)\right)$):
$$\vert G_F(z)-G_{F_0}(z)\vert \le \frac{C_2}{d-1}\Vert F-F_0 \Vert^{\alpha};\;
\forall z\in ({\bf C}^{k+1}\setminus\{0\}).$$\qed\\

The Green function $G_F$ induces a continuous, $\omega$-$p.s.h$ function $g_F$ on $\Hd(\PP)\times \PP$
which will also be called a {\it Green function} of $F$:
\begin{equation}\label{2}  
g_F\circ
  \pi :=G_F - \log \Vert\;\Vert.
\end{equation}

\begin{rem}\label{1aa} {\em It is straightforward to check that   $g_F(\pi(z)) \le \frac{M}{d-1}$, where 
\linebreak $M:=\sup_{||z||=1}||F(z)||$.
In particular} for every compact subset ${\cal K} \subset {\bf C}^{N+1}=\Hd({\bf C}^{k+1})\cup \Sigma_d$, $g_F$ is bounded from 
above on 
$\big({\cal K}\cap\Hd({\bf C}^{k+1})\big)\times \PP$.
\end{rem}
\subsection{Green currents and measures}

Let $f\in \Hd(\PP)$ and $F\in \Hd({\bf C}^{k+1})$   be   a lift of $f$.
One defines a closed, positive $(1,1)$-current $T_f$ on $\PP$ by setting:

\begin{equation}\label{4}
T_f:=dd^c g_F + \omega.
\end{equation}

As \eqref{2} shows, this current may equivalently be defined by $\pi^* T_f = dd^c G_F$.
Since $g_{aF}=\frac{1}{d-1}\log\vert a\vert +g_F$, this current does not depend on the choice
of the lift $F$ and will be called the {\it Green current} of $f$. The functional equation \eqref{3}
implies that:

\begin{equation}\label{5}
f^*T_f=dT_f.
\end{equation}
 
The {Green measure} $\mu_f$ of $f$ is defined by 
$$\mu_f:=\left(T_f\right)^k.$$
It is a
probability measure with respect to which $f$ has constant Jacobian: $f^*\mu_f=d^k\mu_f$.
It follows that $\mu_f$ is $f$-invariant ($f_*\mu_f=\mu_f$)
 and $f$-ergodic.\\

It will also be useful to consider the probability measures $m$ and $\mu_F$ defined on ${\bf C}^{k+1}$ by:
$$m:=\left(dd^c\log^+\Vert\;\Vert\right)^{k+1}\;\;\;\;\;\;\;\;\;\;\;\;\;\mu_F:=\left(dd^c G^+_F
\right)^{k+1};$$
these measures are respectively supported by $\{\Vert\;\Vert=1\}$ and $\{G_F=0\}$; they are related to $\omega^k$ and $\mu_f$
by:
$$\pi_*m=\omega^k\;\;\;\;\;\;\;\;\;\;\;\;\;\pi_*\mu_F=\mu_f.$$

\subsection{ Lyapunov exponents}\label{ssle}

Let $f\in \Hd(\PP)$ and $F\in \Hd({\bf C}^{k+1})$ a lift of $f$.
We shall denote by $L(F)$ the sum of  Lyapunov exponents of $F$ with respect to $\mu_F$ and by
$L(f)$ the sum of  Lyapunov exponents of $f$ with respect to $\mu_f$. The following facts hold:
\begin{itemize}\item[(i)] $L(F)=\int_{{\bf C}^{k+1}} \log \vert det F'\vert \mu_F$;

\item[(ii)] $L(F)=L(f)+\log d$ (see \cite{Jon});

\item[(iii)]$L(F^n)=nL(F)$, for all $n \in {\bf N}^*$ (use (i) and $f_*\mu_f=\mu_f$);

\item[(iv)] $L(f)$ is {\em p.s.h.}  on $\Hd  ({\bf C}^{k+1}) $, as it has been proved in the larger setting of polynomial like mappings 
(see
\cite{DiSi}).

\end{itemize}

\subsection{Green metric on ${\cal O}_{\PP}(D)$}\label{ssGm}

Let $D\in {\bf N}^*$. The line bundle ${\cal O}_{\PP}(D)$
over $\PP$ is conveniently seen as the quotient of $({\bf C}^{k+1}\setminus\{0\}) \times{\bf C}$ by the relation
$(z,x)\equiv(uz,u^Dx)$ for all $u\in {\bf C}^*$,   denoting its elements by   $[z,x]$. The {\it canonical metric} on
 ${\cal O}_{\PP}(D)$ may be expressed by:
$$\Vert[z,x]\Vert_0:=e^{-D\log \Vert z\Vert}\vert x\vert.$$
The homogeneity of $G_F$ allows us to associate to any $F\in \Hd({\bf C}^{k+1})$ a {\it Green metric} defined 
on ${\cal O}_{\PP}(D)$
by:
$$\Vert[z,x]\Vert_{G_F}:
=e^{-DG_F  (z)  
}\vert x\vert.$$
The main interest of endowing ${\cal O}_{\PP}(D)$ with such a metric is to produce a very
 useful formula for $L(F)$. To this purpose we will denote   by   $J_F$ the holomorphic section induced by $det F'$ on 
${\cal O}_{\PP}(D)$ for $D=(k+1)(d-1)$:
$$J_F\circ \pi = [z,det F'(z)];\;\;\;\forall z\in {\bf C}^{k+1}.$$
Then we have the following lemma:
\begin{lem}\label{6}
Let $F\in
\Hd(\PP)$ and $D=(k+1)(d-1)$. Let us endow ${\cal O}_{\PP}(D)$ with the canonical and the Green metrics. Then
the following identities occur:
\begin{itemize}
\item[1)]$L(F)=\int_{{\bf C}^{k+1}} \log \vert det F'\vert\; \mu_F
=\int_{\PP}\log \Vert J_F\Vert_{G_F}\;
 \mu_f$
\item[2)]$
\int_{{\bf C}^{k+1}} \log \vert det F'\vert\;m=\int_{\PP}\log \Vert J_F\Vert_0 \;\omega^k.$
\end{itemize}
\end{lem}

\proof 1) 
 As
$G_F$ identically vanishes on the support of $\mu_F$, we have $\int_{{\bf C}^{k+1}} \log \vert det F'\vert\; \mu_F=
\int_{{\bf C}^{k+1}} \log\left(e^{-DG_F}\vert det F'\vert \right)\mu_F=
\int_{{\bf C}^{k+1}} \log \Vert J_F\circ \pi \Vert \mu_F $ and the conclusion follows from 
 $\pi_*\mu_F=\mu_f$.\\
2) We proceed   in the same way,   using the fact that $\log^+\Vert\;\Vert$ identically vanishes on the support of $m$ and
 $\pi_*m=\omega^k$.\qed\\

\section{A current detecting the holomorphic motion of repulsive cycles}\label{SMvH}

Let $\fL$ be a holomorphic family of endomorphisms of $\PP$ parametrized by a complex manifold $X$. 
The $p.s.h.$ function $L(\lambda)=L(f_{\lambda})$
given by the sum
of       Lyapunov exponents of $f_{\lambda}$,   provides   a closed, positive $(1,1)$-current $dd^c L$ on $X$. In this section, we will
 show that $dd^c L$ vanishes   if   the repulsive cycles of $f_{\lambda}$ move holomorphically. Let us precisely state what we
 mean by this holomorphic motion.
\begin{defi}\label{0a}
The repulsive cycles of period $n$ of $\fL$ move holomorphically over an open subset $U$ of $X$ if and only if there
 exists a collection of holomorphic mappings $\alpha_{n,j}:U\to \PP$ such that, for any $\lambda \in U$, the set of $n$-repulsive
 cycles is given by $\{\alpha_{n,j}(\lambda)\}$.
\end{defi}

In dimension $k=1$, it is well known that the Julia set of $f_{\lambda}$ depends continuously on $\lambda\in U$ if and only if the repulsive cycles of 
sufficiently high
 order of $f_{\lambda}$ move holomorphically on $U$ (see \cite{MacMu}, Theorem 4.2).\\
Although such a phenomenon is far from being
clear in higher dimension, we would like to motivate the study of $dd^c L$ as a bifurcation current by
the following result:

\begin{theo}\label{a}  Let $\left(f_{\la}\right)_{\la \in X}$ be a
holomorphic family of 
endomorphisms of $\PP$ with algebraic degree $d$.
If all repulsive 
cycles of $f_{\la}$ of period $n\ge n_0$ move  holomorphically  on some open subset
$U$ of $X$ then the sum $L(f_{\la})$ of  Lyapunov exponents of $f_{\la}$ is pluriharmonic 
on $U$.
\end{theo}

\proof  We may assume that $U$ is a small ball on which  $\{ f_{\la} \}_{\la\in U}$ lifts to some 
holomorphic family
$\{F_{\la}\}_{\la \in U}$. Then it is not hard to see that
the set of $n$-repulsive cycles of 
$F_{\la}$ 
is given by  
$\{a_{n,j}(\la); \la \in U\}$ 
where the maps $a_{n,j}:U\to {\bf C}^{k+1}$ are holomorphic. The 
number of 
elements of 
$\{a_{n,j}(\la); \la \in U\}$ does not depend on $\la$ and will be 
denoted   by   $N_n$.\\
By a theorem of Briend-Duval (see \cite{BD}, Theorem 2),
the Green measure 
$\mu_{F_{\la}}$ of $F_{\la}$ is the weak limit of a sequence of discrete 
measures: 
 
$$
 \frac{1}{N_n}\Sigma_j\delta_{a_{n,j}(\la)}=:\mu_{F_{\la},n}\to\mu_{F_{\la}}.$$
It is therefore natural to consider the sequence 
of pluriharmonic functions

$$L_n(\la):=\int_{{\bf C}^{k+1}}\log \vert det F'_{\la}\vert
\;\mu_{F_{\la},n}=\frac{1}{N_n}\Sigma_j \log \vert det F'_{\la}\left(a_
{n,j}(\la)\right)\vert$$
and to compare it with $L(\la):=\int_{{\bf C}^{k+1}}\log \vert det F'_{\la}\vert
\;\mu_{F_{\la}}=L\left(f_{\la}\right)+\log d$.\\
However, as the function $\log \vert det F'_{\la}\vert$ is unbounded, this 
comparison is not immediate. As we shall see, the fact that the measures 
$\pi_*\mu_{F_{\la}}$ have local $\alpha$-H\" older potentials is essential
to overcome this difficulty.\\

Let us now enter into details and,  
to this purpose, fix a few notations.
The Green function of $F_{\la}$
will be denoted   by   $G_{\la}$ and for any $\varepsilon>0$ we set
$$W_{\la,\varepsilon}:= 
\{G_{\la}=0\}\cap \{\vert det F'_{\la}\vert \le \varepsilon\}.$$  
We shall   call   $d_{n,j}(\la)$ the holomorphic function $det F'_{\la}
\left(a_{n,j}(\la)\right)$ and introduce the following sequence of 
discrete measures:   
 \begin{equation*}
{\cal L}_{n,\la}:=
 \frac{1}{N_n}\sum_j \log \vert d_{n,j}(\la)\vert
 \delta_{a_{n,j}(\la)}.\end{equation*}
Since $G_{\la}\circ F_{\la}^n=d^n G_{\la}$, the $ F_{\la}^n$-fixed points
$a_{n,j}(\la)$ belong to $\{G_{\la}=0\}$ and thus, according to our 
notations, we have 
 \begin{equation}\label{b} 
L_n(\la)=
 {\cal L}_{n,\la}({\bf C}^{k+1})={\cal L}_{n,\la}
\left({W_{\la,\varepsilon}}\right)+{\cal L}_{n,\la}
\left({W^c_{\la,\varepsilon}}\right).
\end{equation}

We now fix $\la_0\in 
U$, a unit vector
$z_0 \in {\bf C}^N\setminus \{0\}$ and $\rho>0$ such that 
$u_{\theta}:=\la_0+\rho e^{i\theta}z_0$ belongs to U for every $\theta \in
{\bf R}$. Since the functions $L_n(\la)$ are pluriharmonic on $U$, the
identity \eqref{b} may be rewritten as
 \begin{equation}\label{c}
L_n(\la_0)=\frac{1}{2\pi}\int_0^{2\pi}
 {\cal L}_{n,u_{\theta}}
 \left({W^c_{u_{\theta},\varepsilon}}\right)d\theta
+\frac{1}{2\pi}\int_0^{2\pi}
 {\cal L}_{n,u_{\theta}}
 \left({W_{u_{\theta},\varepsilon}}\right)d\theta.
\end{equation}
 
  Since the   function $L$   is   $p.s.h.$ (see Subsection \ref{ssle})  ,   we simply have to 
deduce 
from \eqref{c} that $
\frac{1}{2\pi}\int_0^{2\pi}L(u_{\theta})d\theta \leq L(\la_0)$.
This will require the following lemmas.

\begin{lem}\label{d}
There exists an universal function $M: ]0,\varepsilon_0]\to ]0,1]$ which 
tends 
to $0$ at $0$ and such that
$\frac{1}{2\pi}\int_0^{2\pi}
 {\cal L}_{n,u_{\theta}}
 \left({W_{u_{\theta},\varepsilon}}\right)d\theta \ge {\cal L}_{n,\la_0}
 \left({W_{\la_0,M(\varepsilon)}}\right)$  
 for every $n\in{\bf N}$ and every $0<\varepsilon\le \varepsilon_0$.
\end{lem}

\begin{lem}\label{e}
$\lim_{\varepsilon\to 0}\left(\liminf_n  {\cal L}_{n,\la}
 \left({W^c_{\la,\varepsilon}}\right)\right)=L(\la)$ for every $\la \in X$.
\end{lem}

Using Lemma \ref{d} and the identities \eqref{b}, \eqref{c} we get 
 \begin{eqnarray*} \frac{1}{2\pi}\int_0^{2\pi}
 {\cal L}_{n,u_{\theta}}
 \left({W^c_{u_{\theta},\varepsilon}}\right)d\theta
=L_n(\la_0)-\frac{1}{2\pi}\int_0^{2\pi}
 {\cal L}_{n,u_{\theta}}
 \left({W_{u_{\theta},\varepsilon}}\right)d\theta\\
 \le L_n(\la_0) - 
 {\cal L}_{n,\la_{0}}
 \left({W_{\la_{0},M(\varepsilon)}}\right)={\cal L}_{n,\la_{0}}
 \left({W^c_{\la_{0},M(\varepsilon)}}\right)
\end{eqnarray*} 
then, by Fatou's theorem we have 
$$\frac{1}{2\pi}\int_0^{2\pi}\liminf_n
 {\cal L}_{n,u_{\theta}}
 \left({W^c_{u_{\theta},\varepsilon}}\right)d\theta 
\le
\liminf_n  {\cal L}_{n,\la_0}
 \left({W^c_{\la_0,M(\varepsilon)}}\right)
.$$
Thus, as $\lim_{\varepsilon\to 0} M(\varepsilon)=0$, the inequality $
\frac{1}{2\pi}\int_0^{2\pi}L(u_{\theta})d\theta \le L(\la_0)$ immediately follows 
from 
Lemma \ref{e} when $\varepsilon\to 0$. This ends the proof of Theorem \ref{a}.
\qed\\

\noindent\underline{Proof of lemma \ref{d}}: we shall use the following fact 
which 
is a direct consequence of  Montel's theorem and  Hurwitz lemma.\\

\noindent{\bf Fact}:\ {\it Let $0<\rho<r<R$   and   ${\cal S}_{\varepsilon}:=
\{ \varphi \in {\cal O}(\Delta_r,\Delta_R^*); \inf_{\vert
z\vert=\rho}\vert \varphi(z) \vert=\varepsilon \}$. Let $M(\varepsilon):=\sup_
{\varphi\in 
\cal{S}_{\varepsilon}}
\sup_{\vert z\vert \le \rho}\vert \varphi(z) \vert$. Then $\lim_{\varepsilon\to 
0} M(\varepsilon)=0$ 
and in particular $M(\varepsilon)
\le 1$ for $0<\varepsilon \le \varepsilon_0$.}\\

\noindent Let us observe that the functions $d_{n,j}(\la)$ are uniformly 
locally bounded. This follows 
from the continuity of $G_{\la}(z)$ 
and the previous observation 
that $\{a_{n,j}(\la)\}\subset \{G_{\la}= 0\}$. According to our notations we 
have

 \begin{equation}\label{f}\frac{1}{2\pi}\int_0^{2\pi}d\theta\int_{W_{u_
{\theta},\varepsilon}}
 {\cal L}_{n,u_{\theta}}= \frac{1}{N_n}\sum'_j\frac{1}{2\pi}\int_0^{2\pi}\log 
\vert d_{n,j}(u_{\theta}) 
\vert
{\bf 1}_{\{\vert d_{n,j}\vert \le \varepsilon\}}(u_{\theta})d\theta\end
{equation}

\noindent where $\Sigma'_j$ indicates that we only consider the terms for 
which 
$\inf_{\theta}\vert d_{n,j}(u_{\theta}) \vert\le \varepsilon$. 
By the Fact, all these terms satisfy  
$\vert d_{n,j}(u_{\theta}) \vert\le M(\varepsilon)\le 1$ for $\varepsilon\le 
\varepsilon_0$ 
and $ \vert u_{\theta} - \la_0\vert \le \rho$.
In particular, 
$\frac{1}{2\pi}\int_0^{2\pi}\log \vert d_{n,j}(u_{\theta}) \vert
{\bf 1}_{\{\vert d_{n,j}\vert \le \varepsilon\}}(u_{\theta})d\theta \ge
\frac{1}{2\pi}\int_0^{2\pi}\log \vert d_{n,j}(u_{\theta}) \vert
d\theta=\log \vert d_{n,j}(\la_0) \vert$. Thus \eqref{f} yields 

$$ \frac{1}{2\pi}\int_0^{2\pi}
 {\cal L}_{n,u_{\theta}}
 \left({W_{u_{\theta},\varepsilon}}\right)d\theta \ge \frac{1}{N_n}\sum'_j\log 
\vert 
d_{n,j}(\la_0)\vert.$$
Finally, as $\vert d_{n,j}(\la_0) \vert\le M(\varepsilon)\le 1$ for all terms 
in $\Sigma'_j$, we have   
$$\frac{1}{N_n}\sum'_j\log \vert d_{n,j}(\la_0)\vert \ge 
\frac{1}{N_n}\sum_j \log \vert d_{n,j}(\la_{0}) \vert
{\bf 1}_{\{\vert d_{n,j}\vert \le M(\varepsilon)\}}(\la_{0})=
{\cal L}_{n,\la_{0}}
 \left({W_{\la_{0},M(\varepsilon)}}\right)$$ and the conclusion follows.\qed\\

\noindent\underline{Proof of lemma \ref{e}}: Let us denote   by   $\log_{\varepsilon}x$ 
  an increasing   smooth function on $[0,+\infty[$ such that 
$\log_{\varepsilon}x\ge 2\log\varepsilon$ for $ 0\le x < \varepsilon<1$ and  
$\log_{\varepsilon}x= \log x $
for $x\ge \varepsilon.$ Then

\begin{equation}\label{g}0\le{\cal L}_{n,\la}
 \left({W^c_{\la,\varepsilon}}\right)-\int_{{\bf C}^{k+1}} \log_{\varepsilon}
\vert det F'_{\la}\vert 
\mu_{F_{\la},n}
\le -2(\log \varepsilon) \mu_{F_{\la},n}\left(\{\vert det F'_{\la}\vert\le 
\varepsilon\}\right).\end{equation}
There are constants $a,A>0$ such that $\{\vert det F'_{\la}\vert\le 
\varepsilon\}\subset \pi^{-1}
(V_{A\varepsilon^{a}}(C_{f_{\la}}))$ where $ V_{t}(C_{f_{\la}})$ denotes a $t$-neighbourhood of $
C_{f_{\la}}$. Since $\mu_{\la}:=\pi_{\star}\mu_{F_{\la}}$ has (local) $\alpha$-H\" older continuous 
potential we have $\mu_{\la}
\left(V_{A\varepsilon^{a}}(C_{f_{\la}})\right)\le cst\; \varepsilon ^{a\alpha}
$ (see the proof of Theorem 1.7.3 in \cite{Sb}).
Thus, for $n$ big enough, 

\begin{equation}\label{h}\mu_{F_{\la},n}\left(\{\vert det F'_{\la}\vert\le 
\varepsilon\}\right)
\le 2 \mu_{F_{\la}}\left(\{\vert det F'_{\la}\vert\le 2\varepsilon\}
\right)      
\le cst\; \varepsilon ^{a\alpha}.\end{equation}

From \eqref{g} and \eqref{h} we get 

$$0\le\liminf_n {\cal L}_{n,\la}
 \left({W^c_{\la,\varepsilon}}\right)-\int_{{\bf C}^{k+1}} \log_{\varepsilon}
\vert det F'_{\la}\vert 
\mu_{F_{\la}}
\le -cst\; \varepsilon^{a\alpha}\log\varepsilon$$
and the conclusion follows by making $\varepsilon\to 0$.\qed

  \section {Formulas for the  Lyapunov exponent of a rational function}\label{SF1}
  
 In this section we establish some formulas which relate the Lyapunov exponent to the critical points of a rational function.
 
  \begin{theo}\label{i}
 Let $f$ be a rational function of degree $d$ and $F$ be one of its lifts to
 ${\bf C}^2$. The  Lyapunov exponent $L(f)$ of $f$ is given by one of the 
following 
 formulas:
\begin{itemize}
 \item[(i)]$ L(f)+\log d=\int_{\pp} g_F [C_f] -2(d-1)\int_{\pp} g_F \left
(\mu_f + \omega\right)
+ \int_{\pp}\log 
\Vert {J_F}\Vert_0 \omega.$ 
  \item[(ii)]$L(f)+\log d = \int_{\pp}g_F[C_f]\;-2(d-1)\int_{\pp}g_F\left
(\mu_f+\omega\right)\;+
  \int_{{\bf C}^2}\log\;\vert detF'\vert m.$
 \item[(iii)]If ${\tilde c}_1$,...,${\tilde c}_{2d-2}$ are chosen such that
  $detF'(z)=\Pi_{j=1}^{2d-2}
  {\tilde c}_j\wedge z$ one has:
\item[]$L(f) + \log d = \Sigma_j G_F\left({\tilde c}_j\right)\;
 -(d-1)\big(1 +
  2\int_{\pp}
  g_F\left(\mu_f+\omega\right)\big)
 .$\end{itemize} \end{theo}
  
\proof Let us start with the first formula.
 We know that $L(f)+\log d =L(F)$. We shall use the formalism introduced in the Subsection \ref{ssGm}.
 By the first assertion
 of Lemma
\ref{6}
and the definition of $\mu_f$ we have

\begin{equation}\label{j}L(F) = \int_{\pp}\log 
  \Vert {J_F}\Vert_{G_F} \mu_f = \int_{\pp}\log 
  \Vert {J_F}\Vert_{G_F} dd^c g_F + \int_{\pp}\log 
  \Vert {J_F}\Vert_{G_F} \omega.\end{equation}
  After an integration by parts (see next section for a careful
 justification) the identity \eqref{j} yields 
\begin{equation}\label{k}L(F)=\int_{\pp} g_F dd^c \log 
  \Vert {J_F}\Vert_{G_F} + \int_{\pp}\log 
  \Vert {J_F}\Vert_{G_F} \omega.\end{equation}
  Using   the   Poincar\'e-Lelong equation $dd^c \log\Vert
  {J_F}\Vert_{G_F}=-(2d-2)\mu_f + [C_f]$, 
  \eqref{k} becomes:
   
  \begin{equation}\label{l}L(F)=\int_{\pp} g_F [C_f] -
(2d-2)\int_{\pp} g_F \mu_f +
  \int_{\pp}\log 
  \Vert {J_F}\Vert_{G_F} \omega.\end{equation}
 After observing that $\Vert\cdot\Vert_{G_F}=e^{-(2d-2)g_F}\Vert
  \cdot\Vert_0$ we may rewrite the last integral in \eqref{l} as:
   $\int_{\pp}\log 
  \Vert {J_F}\Vert_0 \omega -(2d-2)\int_{\pp} g_F 
  \omega$ and this gives our first formula.\\

In order to establish the second formula, we simply transform the first one by 
using 
the second assertion of 
 Lemma \ref{6}.\\

Let us finally prove the third formula. Picking $U_j \in {\bf U}(2,{\bf C})$ 
such that
 $U_j^{-1}({\tilde c}_j)=\left(\Vert {\tilde c}_j\Vert,0 \right)$ we have 
$U_j z\wedge {\tilde c}_j=-z_2\Vert {\tilde c}_j\Vert$. Since
$\int_{{\bf C}^2}\log\vert z_2\vert m =-\frac{1}{2}$, we get 
$\int_{{\bf C}^2}\log \vert det F'\vert m = \Sigma_j \int_{{\bf C}^2}\log 
\vert  
U_j z\wedge {\tilde c}_j\vert m = \Sigma_j\log\Vert{\tilde c}_j\Vert - (d-1)
$.\\
On the other hand, $\int_{\pp} g_F [C_f]=\Sigma_j g_F\circ\pi ({\tilde c}_j)= 
\Sigma_j G_F ({\tilde c}_j)-\Sigma_j \log\Vert ({\tilde c}_j)\Vert$. It then 
suffices to 
  replace   these
identities
in the second formula.\qed 
 
\bigskip

 The usefulness of our formulas   consists in   the fact that the function 
\linebreak $B(F):= \int_{\pp} g_F \left(\mu_f +\omega\right) $ 
 is pluriharmonic. This important property is easy to check by considering
 the formulas for both $f$ and $f^2$. Indeed, since $L(f^2)=2L(f)$, $G_{F^2}
=G_F$
  and (consequently) $B(F^2)=B(F)$ one immediately
 obtains 
 a pluriharmonic expression for $B(F)$ by comparison. 
 
\begin{theo}\label{m}
  The function $B(F):=\int_{\pp} g_F \left(\mu_f +
 \omega\right)
 $ is
  pluriharmonic on 
  ${\cal H}_d\left({\bf C}^2\right)$.
 \end{theo}
  
\bigskip

Using the third formula of  Theorem \ref{i} and Theorem \ref{m} we   get   the 
following corollary,   
previously obtained by DeMarco
(see \cite{Demar}). It allows to relate the pluriharmonicity of $L
(F)$ with the stability of the 
dynamic of critical points.
As we shall see in Section \ref{sbc}, this is a key point when approaching the Ma\~n\'e-Sad-Sullivan theory via 
potential-theoretic methods.

\begin{cor}\label{m1} Let $\fL$ be a holomorphic family of rational maps of 
degree $d$.  Then $dd^cL(f_{\la})=dd^c\sum_{j=1}^
{2d-2}G_F(\tilde 
c_j)$.
\end{cor}
 
 Using Proposition \ref{1a} (ii),
 one also reads on the third formula of Theorem \ref{i} that the  Lyapunov 
exponent
 $L(F)$ is an H\" older-continuous function in $F$
. The continuity was first proved by Ma\~n\'e \cite{Mane}.    

 \begin{cor}\label{n} The  function $L(F)$ is $p.s.h.$ and H\" older-continuous 
on ${\cal
 H}_d\left({\bf C}^2\right)$.
  \end{cor}
 
 \noindent\underline{Proof of theorem \ref{m}}.
 We may consider a local holomorphic parametrization $\la\mapsto
 F_{\la}$ of ${\cal H}_d\left({\bf C}^2\right)$
 defined on some open subset $U$ of 
 ${\bf C}^{2d+2}$. We shall denote   by  
 $f_{\la}$ the induced map on ${\bf P}^1$ and set
 $B(\la):=B\left(F_{\la}\right)$.
 There exists an analytic subset $A$ of $U$ such that, for any $\lambda\in U\setminus A$, the critical points of
 $f_{\la}$ consist in $2d-2$ distinct, regular values of $f_{\la}$.
 As the function $B(F)$ is locally bounded, it suffices to show that 
it is 
 pluriharmonic on any sufficiently small ball contained in $U\setminus A$.\\
  On such a ball $B$, there are $2d-2$ holomorphic maps $\tilde{c}_j$ such
that
    $$detF'_{\la}=\Pi_{j=1}^{2d-2}\;\;{\tilde c}_j(\la)\wedge z.$$
  Moreover, for each $1\le j \le 2d-2$, there are $d$ holomorphic maps 
  ${\tilde c}_{j,i}$ 
such that $F_{\la}\circ{\tilde c}_{j,i}(\la)={\tilde c}_j(\la)$ 
and therefore:
 $$detF^{2'}_{\la}=h(\la)
  \left(\Pi_{j=1}^{2d-2}\;\Pi_{i=1}^{d}\;\;{\tilde c}_{i,j}(\la)\wedge
  z\right)
  \left(\Pi_{j=1}^{2d-2}\;\;{\tilde c}_j(\la)\wedge z\right)$$
where $h$ is a non-vanishing holomorphic function on $B$. Let $N$ denote the 
degree of
$det F^{2'}_{\la}$, after setting
 ${\tilde c}'_{j,i}=h^{\frac{1}{N}}\tilde{c}_{j,i}$ and 
${\tilde c}'_{j}=h^{\frac{1}{N}}\tilde{c}_{j}$ we get
 $$detF^{2'}_{\la}= \left(\Pi_{j=1}^{2d-2}\;\Pi_{i=1}^{d}\;\;{\tilde
 c}'_{i,j}(\la)\wedge z\right)
 \left(\Pi_{j=1}^{2d-2}\;\;{\tilde c}'_j(\la)\wedge z\right).$$

We are now in order to use the third formula of Theorem \ref{i} for $f_{\la}^2$.
  Since   $G_{F^2}=G_F$, $\mu_{f^2}=\mu_f$ and
 $B\left(F^2\right)=
 B\left(F \right)$,   it yields:

\begin{eqnarray*}
 L\left(f^2_{\la}\right)+\log d^2 +(d^2-1)\left(2B\left(F_\la\right) +1
\right) =
\Sigma_{j=1}^{2d-2}\Sigma_{i=1}^{d}
 G_F\left(\tilde{c}'_{j,i}(\la)\right)
 + \Sigma_{j=1}^{2d-2}G_F\left(\tilde{c}'_{j}(\la)\right)\\
= \log\vert h(\la)\vert  
 + \Sigma_{j=1}^{2d-2}\Sigma_{i=1}^{d}
 G_F\left(\tilde{c}_{j,i}(\la)\right)
 + \Sigma_{j=1}^{2d-2}G_F\left(\tilde{c}_{j}(\la)\right)\\
=\log\vert h(\la)\vert  
 + \Sigma_{j=1}^{2d-2}\frac{1}{d}\Sigma_{i=1}^{d}
 G_F\circ F\left(\tilde{c}_{j,i}(\la)\right)
 + \Sigma_{j=1}^{2d-2}G_F\left(\tilde{c}_{j}(\la)\right)\\
=\log\vert h(\la)\vert  
  + 2\Sigma_{j=1}^{2d-2}G_F\left(\tilde{c}_{j}(\la)\right).
\end{eqnarray*}
 
On the other hand, for $f_{\la}$, the same formula gives:
 $$L\left(f^2_{\la}\right)+\log d^2
=2\left(L\left(f_{\la}\right)+\log d\right)=  
 2\Sigma_{j=1}^{2d-2}G_F\left(\tilde{c}_{j}(\la)\right)-2(d-1)
\left(2B\left(F_\la\right) +1\right).$$
 
 By comparison we thus obtain 
 $2B(\la)+1=\frac{1}{(d-1)^2}\log\vert h(\la)\vert$.\qed\\
 
\begin{rem}\label{n1} {\em By using its pluriharmonicity, one may show that the function
 $B(F)$ is given by
$B(F)=\frac{1}{d(d-1)}\log\vert Res(F)\vert - \frac{1}{2}$. This gives again DeMarco's formula   (\cite{Demar} Corollary 1.6)
  and will be proved in Proposition \ref{z} in arbitrary dimension.}

\end{rem}

\section{A formula for the sum of  Lyapunov exponents of holomorphic 
endomorphisms of 
$\PP$}\label{SF2}

Our aim here is to generalize the results of the previous section to endomorphisms 
of 
$\PP$. We first establish a formula which relates the sum of   the    Lyapunov 
exponents $L(f)$ with the Green current and the current of integration on the 
critical set. 
This extends Theorem \ref{i} (i). We then generalize Theorem \ref{m} and, in particular, obtain
an intrinsic expression for $dd^c L$.
\medskip

\begin{theo}\label{o}
Let $f$ be a holomorphic endomorphism of $\PP$   of   algebraic degree 
 $d \geq 2$. Let $F$ be one of the lifts
of $f$ to ${\bf C}^{k+1}$ and $T_f=dd^cg_F+\omega$ be the Green current of 
$f$. 
Then the sum of   the    Lyapunov exponents $L(f)$
of $f$ is given by:
\begin{equation}\label{o1}L(f)+\log d=L(F)=H(F)-(k+1)(d-1)B(F)\end{equation}
where
\begin{equation}\label{o2}H(F):=\sum_{j=0}^{k-1}  \int_{\PP}g_F \ T_f^j \wedge \omega^{k-j-1}\wedge [C_f] +\int_{{\bf P}^k}\log||J_F||_0\omega^k
\end{equation}
and
\begin{equation}\label{o3}B(F):=\sum_{j=0}^k\int_{{\bf P}^k}g_F \  T_f^j \wedge \omega^{k-j}.\end{equation} 
\end{theo}
\proof According to Lemma \ref{6} we have
$$L(f)+\log d=\int_{{\bf P}^k} \log||J_F||_G \ T_f^k.$$
Let us start by showing that
$$L(f)+\log d=\sum_{j=0}^{k-1}\int_{{\bf P}^k}\log||J_F||_G\ dd^cg_F \wedge 
T_f^j \wedge 
\omega^{k-j-1}+$$
\begin{equation}\label{p}+\int_{{\bf P}^k}\log||J_F||_G\ \omega^k.\end{equation}
  To this purpose   we first note that each term in the above sum is finite (this 
follows immediately 
from the 
Chern-Levine-Nirenberg inequalities) and then   we   observe that:
$$dd^cg_F \wedge \big(\sum_{j=0}^{k-1}T_f^j \wedge \omega^{k-j-1} \big)=
(T_f-\omega)\wedge \big(\sum_{j=0}^{k-1}T_f^j \wedge \omega^{k-j-1} \big)=$$
\begin{equation}\label{p1}=\sum_{j=0}^{k-1}T_f^{j+1} \wedge \omega^{k-j-1} -\sum_{j=0}^{k-1}T_f^j 
\wedge \omega^{k-
j} =T_f^k-\omega^k.\end{equation}
We shall now use the following integration by part 
property which will be 
  proved   separetely.

{\bf Fact:}   for $0 \le j <k$,  
$$\int_{{\bf P}^k}\log||J_F||_G \ dd^cg_F\wedge T_f^j \wedge \omega^{k-j-1}
=\int_{{\bf P}^k}g_F \ dd^c\log||J_F||_G \wedge  T_f^j \wedge \omega^{k-j-
1}.$$

\medskip

This allows us to transform the identity \eqref{p} and get:
$$L(f)+\log d=\sum_{j=0}^{k-1}\int_{{\bf P}^k}g_F \ dd^c\log||J_F||_G \wedge  
T_f^j \wedge 
\omega^{k-j-1}+$$
$$+\int_{{\bf P}^k}\log||J_F||_G\  \omega^k.$$
Next, by the Poincar\' e-Lelong equation $dd^c\log||J_F||_G=[C_f]-(k+1)(d-1)
T_f$, we obtain:
\begin{equation}\label{q}L(f)+\log d=\sum_{j=0}^{k-1}\int_{{\bf P}^k}g_F \ 
T_f^j \wedge 
\omega^{k-j-1}\wedge [C_f]-
\end{equation}
$$-(k+1)(d-1)\sum_{j=0}^{k-1}\int_{{\bf P}^k}g_F \  T_f^{j+1} \wedge \omega^{k-
j-1}
+\int_{{\bf P}^k}\log||J_F||_G\omega^k.$$
Finally, as $||\ ||_G=e^{-(k+1)(d-1)g_F}||\ ||_0$, we may replace the last 
integral in 
\eqref{q} by 
$\int_{{\bf P}^k}\log||J_F||_0\omega^k-(k+1)(d-1)\int_{{\bf P}^k}g_F \ 
\omega^k$ and this 
immediately yields to the expected formula.
\qed\\

It remains to establish the Fact. We shall proceed by regularization and use the following lemma.

\begin{lem}\label{r} Let $\{\phi_n\}_{n \in {\bf N}^*}$ be a decreasing 
sequence of increasing 
smooth convex functions
on ${\bf R}$ such that:
  \begin{itemize} \item[(i)] $\phi_n(x)=-n,\  on \ \ ]-\infty,-n-\frac{1}{n}]$
\item[(ii)]
$\phi_n(x)=x,\  on \ \ ]-n+\frac{1}{n},+\infty[$ 
\end{itemize}

\noindent Let $\log_nx$ be defined by $\log_nx:=\phi_n(\log x)$. 
Then $\{\log_n ||J_F||_0 \}_{n \in {\bf N}^*}$ is a decreasing 
sequence of smooth functions, which converges to 
$\log||J_F||_0$. Moreover \linebreak $dd^c \log_n \Vert J_F\Vert_0 + (k+1)(d-1)\omega\ge 0$ for all $n\in {\bf N}^*$.
\end{lem}

\medskip

\noindent The proof is   a   straightforward computation and we omit it. 

\medskip

\noindent\underline{Proof of the Fact}.    Using   the relation
$||\ ||_{G_F}=e^{-(k+1)(d-1)g_F}||\ ||_0$,   we get: 
$$\int_{{\bf P}^k}\log||J_F||_G \ dd^cg_F\wedge T_f^j \wedge \omega^{k-j-1}=$$
$$=\int_{{\bf P}^k}\log||J_F||_0 \ dd^cg_F \wedge T_f^j \wedge \omega^{k-j-1}
-(k+1)(d-1)\int_{{\bf P}^k}g_F \ dd^cg_F \wedge T_f^k \wedge \omega^{k-j-1}.$$ 
 Lemma \ref{r} allows us to
use monotone
convergence theorem  (\cite{Sb} Theorem A.6.2), thus   
  $$\int_{{\bf P}^k}\log||J_F||_0 \ dd^cg_F \wedge T_f^j \wedge \omega^{k-j-
1}=\lim_{n \to \infty}\int_{{\bf P}^k}\log_n||J_F||_0\ dd^cg_F \wedge 
T_f^j \wedge 
\omega^{k-j-1}=$$
$$=\lim_{n \to \infty}\int_{{\bf P}^k}g_F\  dd^c\log_n||J_F||_0
\wedge T_f^j \wedge 
\omega^{k-j-1}
=\int_{{\bf P}^k}g_F\  dd^c\log||J_F||_0 \wedge T_f^j \wedge 
\omega^{k-j-1}.$$  
\qed
\bigskip

Our aim now is to compute $dd^cL(f_{\la})$ when
 $\fL$ is a holomorphic family of endomorphisms of $\PP$.
We need   the following technical Proposition
which will be   proved   in the Appendix.

\bigskip

\begin{prop}\label{t} Let $X^m \stackrel{\pi}{\longrightarrow} Y^n$ be a holomorphic 
submersion between complex 
manifolds.
If $R$ is a current on $X$, for $y \in Y$   the slice (if it exists) 
of $R$ along the fiber $\pi^{-1}(y)$ is denoted by  $R_y$. Let $u_1, \dots, u_h$ be almost plurisubharmonic, 
locally bounded functions on $X$ and $T$ be a positive, 
closed $(k,k)$-current on $X$,
with $h+k \leq m-n$. Thus, for a.e. $y \in Y$,
$$(u_1 dd^c u_2 \wedge \dots \wedge dd^c u_h \wedge T)_y=$$
$$=u_{1 |\pi^{-1}(y)} dd^c (u_{2|\pi^{-1}(y)}) \wedge \dots \wedge dd^c (u_{h|\pi^{-1}(y)}) \wedge T_y.$$
\end{prop}
Let us recall that, for a $(k,k)$-current $R$ on $X$, slicing is characterized (for a.e. $y \in Y$)  by the following identity:
\begin{equation}\label{u}
\int_X R \wedge \psi \wedge \pi^*\phi=\int_Y \big(\int_{\pi^{-1}(y)}R_y \wedge \iota^*_y \psi \big)\  \phi\end{equation}
for every smooth
$(n,n)$-form $\phi$ on $Y$ and for every smooth and compactly supported 
$(m-n-k,m-n-k)$-form $\psi$ on $X$ (here 
$\iota_y:\pi^{-1}(y) \to X$ is the inclusion.)

\bigskip

By Theorem \ref{o},   $L(f_{\la}) + \log d=H(F_{\la})-(k+1)(d-1)B(F_{\la})$.  We first compute $dd^cH$:

\medskip

\begin{prop}\label{u1} Let $\fL$ be a holomorphic family of endomorphisms of $\PP$ such that there is
a holomorphic lift $\{F_{\la }\}_{\la \in X}$ to ${\bf C}^{k+1}$. Then

$$dd^cH(F_{\la })=p_*\big((dd^cg_{F_{\la }}+\omega)^k \wedge [C_X]\big)$$
where   \label{NB 6} $C_X$ is the hypersurface of $X \times \PP$ defined by the equation {\em det}$F'_{\la}(z)=0$   and 
\linebreak $p:X\times{\bf P}^k \to X$ is the 
canonical projection. \end{prop}

\noindent{\bf Remark.} $dd^cg_{F_{\la}}$ involves derivatives in both $\la \in X$ and  $z \in \PP$.

\medskip

\proof Let $q=dim_{\bf C}X$, for a $(q-1,q-1)$-form $\phi$
with compact support on $X$ we have
$$<dd^cH,\phi>= \int_X \big(\int_{{\bf P}^k}g_{F_{\la}}(\sum_{j=0}^{k-1}(dd^cg_{F_{\la}}+\omega)^j\wedge \omega^{k-j-1}) 
\wedge[C_{f_{\la}}] \big)\  dd^c\phi \ +$$
$$+ \int_X \big(\int_{{\bf P}^k}\log||J_{F_{\la}}||_0 \omega^k \big)\ dd^c  \phi.$$
Since $[C_{f_{\la}}]$ is the slice of $[C_X]$ (see \cite{Siu} (10.4)), by means of the 
Proposition \ref{t} the first integral is
$$\int_{X \times {\bf P}^k}p^*\phi \wedge dd^cg_{F_{\la}} \wedge (\sum_{j=0}^{k-1}(dd^cg_{F_{\la}}+\omega)^j\wedge \omega^{k-j-1}) 
\wedge[C_X].$$
  By   Poincar\' e - Lelong formula 
 $[C_X]=dd^c\log||J_{F_{\la}}||_0+(k+1)(d-1)\omega$ one sees that
 $\log||J_{F_{\la}}||_0$ is almost plurisubharmonic and therefore locally summable. 
Thus the second integral is
$$\int_{X \times {\bf P}^k}p^*\phi \wedge dd^c \log||J_{F_{\la}}||_0 \wedge \omega^k=$$
$$=\int_{X \times {\bf P}^k}p^*\phi \wedge \omega^k \wedge [C_X]-
(k+1)(d-1)\int_{X \times \PP}p^*\phi \wedge \omega^{k+1}.$$
But $\omega^{k+1}=0$ on $X\times {\bf P}^k$ and therefore, after summing up, we obtain:
$$<dd^cH,\phi>=\int_{C_X}p^*\phi \wedge (dd^cg_{F_{\la}}+\omega)^k=<(dd^cg_{F_{\la}}+\omega)^k \wedge[C_X],p^*\phi>.$$
\qed

\bigskip

 Now we can also extend  Theorem \ref{m} to the $k$-dimensional case. We shall use the same
device, that is to compare formulas for $F_{\la}$ and $F_{\la}^2$.
 
\medskip
\begin{theo}\label{s}  The function $B(F)$ is pluriharmonic on
 ${\cal H}_d({\bf C}^{k+1})$.\end{theo}

\proof Let us start with a   claim  :

\medskip 

\noindent{\bf Claim:} {\em $H(F)$ is {\em p.s.h.} on ${\cal H}_d({\bf C}^{k+1})$ and} \ 
$dd^cH(F^2)=2dd^cH(F).$

\noindent\underline{Proof of the Claim}. Let $X:={\cal H}_d({\bf C}^{k+1})$, the projection $X \ni F \mapsto f \in {\cal H}_d({\bf P}^k)$
 defines a holomorphic family $\{f\}_{F \in X}$ and the plurisubharmonicity follows from  Proposition \ref{u1}.

   \label{NB 7}Let $C_X$ be as above and denote by $C'_X$ the analogous {\em critical set} of the family $\{f^2\}_{F \in X}$. 
Considering the map $ \Phi:X\times {\bf P}^k \to X\times {\bf P}^k$
defined by $\Phi(F, z):=(F,f(z))$, we get $[C'_X]=[C_X]+\Phi^*[C_X]$.

As $F^2$ and 
$F$ have the same Green function $g_F$, Proposition \ref{u1} gives
$$dd^cH(F^2)=p_*\left((dd^cg_F+\omega)^k \wedge [C_X']\right).$$
From $G_F\circ F= d.G_F$, it follows $\Phi^*(dd^cg_F+\omega)=d.(dd^cg_F+\omega)$, thus
$$dd^cH(F^2)=p_*\left((dd^cg_F+\omega)^k \wedge [C_X]\right)+
\frac{1}{d^k}p_*\Phi^*\left((dd^cg_F+\omega)^k \wedge [C_X]\right).$$
But $p \circ \Phi=p$, thus $p_*=p_*\Phi_*$; moreover $\Phi_*\Phi^*=d^k$id, thus
$$\frac{1}{d^k}p_*\Phi^*\left((dd^cg_F+\omega)^k \wedge [C_X]\right)=p_*\left((dd^cg_F+\omega)^k \wedge [C_X]\right),$$
therefore $dd^cH(F^2)=2dd^cH(F)$. \qed

\bigskip

\noindent\underline{End of the proof of Theorem \ref{s}}. 
Since $L(F^2)=2L(F)$ we have
$$dd^cL(F^2)=2dd^cL(F)=2\big( dd^cH(F)-(k+1)(d-1)dd^cB(F) \big)$$
on the other hand, since $B(F^2)=B(F)$, we may use the Claim and get:
$$dd^cL(F^2)=dd^cH(F^2)-(k+1)(d^2-1)dd^cB(F^2)=$$
$$=2dd^cH(F)-(k+1)(d^2-1)dd^cB(F).$$
By comparison we get $(d-1)^2dd^cB(f)=0$, thus $B$ is pluriharmonic on ${\cal H}_d({\bf C}^{k+1})$. 
\qed
\bigskip

\begin{cor}\label{v} Let $\fL$ be a holomorphic family of
 endomorphisms of ${\bf P}^k$ with algebraic 
degree $d\geq 2$. Then
$$dd^cL(f_{\la})=p_*((dd^cg_{F_{\la}}+\omega)^k\wedge[C_X])$$
and on $X \times {\bf P}^k$
$$(dd^cg_{F_{\la}}+\omega)^{k+1}=0.$$
\end{cor}
\medskip 

\begin{rem}\label{w} {\em As we have already noted, the operator $dd^c$ in the above formula involves derivatives in both $\la \in X$
 and 
$z \in {\bf P}^k$; thus the current $\tilde T:=dd^cg_{F_{\la}}+\omega$ is different from
the Green current.    The current $\tilde T$  
 depends only on the family $\{f_{\la}\}$ and not on the local lift $\{F_{\la}\}$; 
moreover it is positive on $X \times \PP$ since $G_{F_{\la}}(z)$ is {\em p.s.h.} on $X \times  ({\bf C}^{k+1} \setminus \{0\})$
(see Proposition \ref{1a}).
Using    this   current  we may express  the formulas of Corollary \ref{v} in a synthetic way and 
avoid any reference to the lift $\{F_{\la}\}$, which in general is only defined  
locally:
$$dd^cL(f_{\la})=p_*(\tilde T^k\wedge[C_X])$$ and $$\tilde T^{k+1}=0.$$ }\end{rem}

\bigskip

\noindent\underline{Proof of Corollary \ref{v}}.  
From Proposition \ref{u1} and Theorem \ref{s} we get the first statement. 

We argue as in the proof of  Proposition \ref{u1} (using again Proposition \ref{t})  ;   choosing an open subset
$V \subset X$ such that there is a holomorphic family of lifts $\{F_{\la}\}_{\la \in V}$  ,   we have
$$<dd^cB,\phi>= \int_V \Big(\int_{{\bf P}^k}g_{F_{\la}}(\sum_{j=0}^{k}(dd^cg_{F_{\la}}+\omega)^j\wedge \omega^{k-j})\Big)
dd^c\phi=$$
$$=\int_{V\times \PP}p^*(\phi)\wedge dd^cg_{F_{\la}}\wedge
\Big(\sum_{j=0}^{k}(dd^cg_{F_{\la}}+\omega)^j\wedge \omega^{k-j}\Big).$$
Then, as in \eqref{p1} we get
$$<dd^cB,\phi>=\int_{V\times \PP}p^*(\phi)\wedge \left( \tilde T^{k+1}- \omega^{k+1}\right)=<p_*(\tilde T^{k+1}),\phi>,$$
because $\omega^{k+1}$ vanishes. Finally, since $dd^cB=0$ and $\tilde T$ is positive, we get $\tilde T^{k+1}=0.$
\qed

\bigskip

By Theorem \ref{s}, the function $B$ is pluriharmonic: this suggests the existence 
of a simpler analytic expression for $B$,
as indeed  Proposition \ref{z} states. Since $B$ is defined by means of dynamical
quantities,   this result seems of some interest.
We shall need the following lemma.

\begin{lem}\label{y} $H^1({\cal H}_d(\PP);{\bf R})=0.$\end{lem}
\noindent {\bf Proof.}  The Fubini-Study form $\omega$ generates
$H^{2N-2}({\bf P}^N;{\bf R})$ and 
$(\iota^*(\omega^{N-1}),\Sigma_d) =$ \linebreak $=\int_{\Sigma_d}\omega^{N-1} = vol(\Sigma_d) \not=0$, where $\iota: \Sigma_d \to {\bf P}^N$ 
is the inclusion;
therefore the map 
${\bf R}= H^{2N-2}({\bf P}^N;{\bf R}) 
\stackrel{\iota^*}{\to} H^{2N-2}(\Sigma_d;{\bf R})={\bf R}$
is an isomorphism. Hence, from the exact sequence
$$H^{2N-2}({\bf P}^N;{\bf R}) \stackrel{\iota^*}{\to}
 H^{2N-2}(\Sigma_d;{\bf R}) \to H^{2N-1}({\bf P}^N, \Sigma_d; {\bf R})
 \to H^{2N-1}({\bf P}^N;{\bf R})=0, $$
it follows that 
$H^{2N-1}({\bf P}^N, \Sigma_d; {\bf R})=0.$

Observe that $\Sigma_d$ is an euclidean neighbourhood
retract (see \cite{D} Prop. IV.8.2, \linebreak VIII.6.12, VIII.7.2) thus 
$H^j({\bf P}^N, \Sigma_d; {\bf R}) = H_{2N-j}({\bf P}^N \setminus \Sigma_d; {\bf R})$. 
In particular \linebreak
$0=H^{2N-1}({\bf P}^N, \Sigma_d; {\bf R})=H_1({\bf P}^N \setminus \Sigma_d;{\bf R})$
and then also its dual   space   $H^1({\bf P}^N \setminus \Sigma_d;{\bf R})$ vanishes. \qed

\bigskip

Now we can establish:

\begin{prop}\label{z} There exists a constant $C_{d,k}$ such that,
 $\forall F \in{\cal H}_d({\bf C}^{k+1})$,
$$B(F)=\frac{1}{d^k(d-1)}\log|Res(F)|+C_{d,k}.$$\end{prop}
\proof If $a\in {\bf C}\setminus\{0\}$, then $g_{aF}=\frac{1}{d-1}\log |a|+g_F$; moreover since
$\int_{{\bf P}^k}T_f^j\wedge\omega^{k-j}=1$ we have
$$B(aF)=\frac{(k+1)}{d-1}\log|a|+B(F).$$
The polynomial $Res(F)$ is homogeneous of degree $(k+1)d^k$, thus the function \linebreak
$d^k(d-1)B(F)-\log|Res(F)|$ is homogeneous of degree $0$ and defines a 
pluriharmonic function $\Phi:{\cal H}_d({\bf P}^k) \to {\bf R}$ such that
$$\forall F\in{\cal H}_d({\bf C}^{k+1}),\ \  \Phi\circ \pi(F)=d^k(d-1)B(F)-\log|Res(F)|.$$
Let ${\cal PH}$ be the sheaf of pluriharmonic functions, by means of Lemma \ref{y}, from the exact sequence
$0 \to {\bf R} \stackrel{i}{\to} {\cal O} \stackrel{Re}{\to} {\cal PH} \to 0$ we get that
 $H^0({\cal H}_d({\bf P}^k), {\cal O}) \stackrel{Re}{\to} H^0({\cal H}_d({\bf P}^k),{\cal PH})$ is surjective;
therefore there exists a holomorphic function $\varphi$ on ${\cal H}_d({\bf P}^k)$ such that
$Re(\varphi)=\Phi$. Setting
 $\psi: = e^{\varphi}$, we obtain
$$\log|\psi| = \Phi \ \hbox{ on} \  {\cal H}_d({\bf P}^k).$$
Using Remark \ref{1aa} one sees that $B$
 is bounded from above on   ${\cal K} \cap {\cal H}_d({\bf C}^{k+1}),$    for every compact ${\cal K} \subset {\bf C}^{N+1}$. It follows that
$Res(F).\psi(\pi(F))$  is locally bounded and thus can be extended to a holomorphic function $\chi$ on ${\bf C}^{N+1}$.
But $\chi$ is clearly homogeneous with the same degree   as   $Res(F)$, thus $\chi$ is a polynomial on ${\bf C}^{N+1}$ and $\psi$
is a constant. \qed

\bigskip

\begin{prop}\label{aa}    The constant $C_{d,k}$ does not depend on $d$, indeed
 $$C_{d,k}=-\frac{1}{2}(k+\frac{k-1}{2}+ \dots + \frac{2}{k-1}+\frac{1}{k}). $$
\end{prop}
\proof See Appendix.
 \bigskip

\begin{rem}\label{ab}
{\em In the one-dimensional case Propositions \ref{z} and \ref{aa} give a new proof  of DeMarco's formula
(see \cite{Demar} Corollary 1.6) since}
$$\int_{{\bf P}^1}g_F(\mu_f+\omega)=\frac{1}{2}(\log|Res(F)|-1).$$\end{rem}

\bigskip

\section{The bifurcation currents}\label{sbc}

In this section, we associate to any holomorphic family $\fL$ in $\Hd(\pp)$ a 
collection
of bifurcation currents $(dd^cL(f_{\la}))^p$ where $1 \leq p \leq dim_{\bf C}
X$. Our main goal is to 
give a rather
precise description of   their   supports and, more precisely, to compare them with 
the hypersurfaces 
consisting of mappings having neutral cycles.
The extremal cases $p=1$ and $p=2d-2$ are of special interest. For $p=1$, we 
partially recover 
Ma\~n\'e-Sad-Sullivan work. For
$p=2d-2$, our description will become significant in the last section when 
introducing a bifurcation 
measure on the moduli space
${\cal M}_d$. Let us   notice   that we shall proceed by induction on $p$.

In order to state the results we must precise a few notations. 

\begin{defi}\label{za} We will consider a holomorphic family $\fL$ of 
elements of $\Hd(\pp)$ parametrized   by   an arbitrary complex manifold $X$. We 
set $D:=2d-2$ 
and denote by $L(\la)$
the {\em p.s.h.} function on $X$ defined by $L(\la):=L(f_{\la})$. Next we 
introduce the following 
subsets of $X$  :  
\begin{itemize} \item[]
${\cal R}:= \{\la_0 \in X;$  the repulsive cycles of sufficiently high period of
$f_{\la}$ move holomorphically on a  fixed neighbourhood $U_0$ of $\la_0 \},$
\item[]
${\cal S}:=\{\la_0 \in X; \la \to f_{\la}^n(C_{f_{\la}})$ is equicontinuous at 
$\la_0\},$
\item[]
  $Per(X,n,e^{2i\pi \theta}):=\{\la_0 \in X;  f_{\la_0}$ has a cycle of period 
$n$ and multiplier 
$e^{2i\pi \theta}\}$, where $\theta \in ]0,1[$. 
\end{itemize}
\end{defi}

It may happen that $Per(X,n,e^{2i\pi \theta})$ is empty or coincides with $X$; 
  otherwise it is a hypersurface 
of $X$.  The union of the irreducible components of codimension $1$ of $Per(X,n,e^{2i\pi \theta})$ will be denoted
by $Per_1(X,n, e^{2i\pi \theta})$.    
For any dense subset ${ E}$ of $]0,1[$, we set
$${\cal Z}_1(X,E)
= \bigcup_{n\in {\bf N}^*, \theta \in E}Per_1
(X,n,e^{2i\pi \theta})$$
Let us recall that the set ${\cal R}$ has been implicitly considered in 
Theorem \ref{a}  ,   which may be stated as
${\cal R}\cap Supp(dd^cL)=\emptyset$. Note also that, in the definition of 
${\cal S}$, the maps 
$\la \to f_{\la}^n(C_{f_{\la}})$ are   considered   as finitely valued holomorphic  maps 
from $X$ to $\pp$.

\medskip

Our description of $Supp(dd^cL)$ contains a substantial part of Ma\~n\'e-Sad-
Sullivan theory (see 
\cite{MSS}).
The originality here relies on the potential-theoretic nature of our proof.

\begin{theo}\label{zb} Let $E$ be a dense subset of $]0,1[$.
Let $\fL$ be a holomorphic family of rational maps of degree $d$ on $\pp$. 
Then
$$\overline{\cal Z}_1(X,E
)={\cal R}^c=Supp(dd^cL)={\cal S}^c.$$
\end{theo}

Let us briefly sketch the proof before entering into details. The inclusion
${\cal S}^c \subset Supp(dd^cL)$ is a consequence of Corollary \ref{m1} and 
was already observed 
by DeMarco (\cite{Demar},
Theorem 1.1). The inclusion $Supp(dd^cL)\subset {\cal R}^c$ 
was proved in 
Theorem \ref{a}   (we recall that the main
ingredient was the equidistribution of repulsive cycles).   The inclusions ${\cal 
R}^c \subset 
\overline{\cal Z}_1(X,E
) \subset {\cal S}^c$
are classical since Ma\~n\'e-Mad-Sullivan work. Their proofs, which we reproduce 
here for sake of 
completeness, only use
the elementary fact that any attractive basin contains at least a critical 
point.

$\underline{{\cal S}^c \subset Supp(dd^cL)}$: Let $\Omega$ be an open ball in 
$X$ on which $L$ is 
pluriharmonic; we have to show that 
$\Omega \subset {\cal S}$.   Shrinking 
$\Omega$
if necessary, we find a $D$-valued holomorphic map $\la \mapsto \tilde C_{f_{\la}}
$ from $\Omega$ to 
${\bf C}^2 \setminus \{0\}$ such that $\pi \circ \tilde C_{f_{\la}}=C_{f_{\la}}$, and an analytic subset $A$ of $\Omega$ 
such that  $\tilde C_{f_{\la}}=\{\tilde c_1(\la), \dots ,\tilde c_D(\la)
\}$ where the
$\tilde c_j(\la)$ are holomorphic maps on $\Omega \setminus 
A$.  

\noindent The product $\Pi(z \wedge \tilde c_j(\la))$ is a well  defined $D$-homogeneous 
polynomial on ${\bf C}^2$
whose coefficients are bounded holomorphic functions on $\Omega \setminus A$. 
It therefore 
coincides with the restriction of a
polynomial $H$ with holomorphic coefficients on $\Omega$. Moreover, as $H$ is 
obviously 
proportional to $detF'_{\la}$,
there exists a non-vanishing holomorphic function $\varphi$ on $\Omega$ such 
that 
$H=\varphi(\la)detF_{	\la}'$.
Thus, after replacing $\tilde c_j(\la)$ by $\left(\varphi(\la)\right)^{1/D}
\tilde c_j(\la)$, we may assume that
\begin{equation}\label{zc} H=\Pi(z \wedge \tilde c_j(\la))=detF_{\la}';\ 
\forall z \in {\bf C}^2,
\forall \la \in \Omega.\end{equation}
In the same way, we may construct a sequence of $D$-homogeneous polynomials 
$H_n$ of the form
\begin{equation}\label{zd}H_n:=h(\la)^{-d^n}\Pi \big(z \wedge F^n_{\la}(\tilde 
c_j(\la))\big)\end{equation}
where $h$ is a non-vanishing holomorphic function on $\Omega$. We will see 
that for a good choice 
of $h$
the coefficients of $H_n$ are uniformly bounded holomorphic functions on 
$\Omega$. As 
$\pi(\{H_n(\la,.)=0\})=f_{\la}^n(C_{f_{\la}})$, this implies that $\Omega 
\subset {\cal S}$.

\noindent Let us to construct $h$. Since $dd^cL=0$ on $\Omega$, it follows 
from \eqref{zc} 
and Corollary \ref{m1} 
that the function $\sum_{j=1}^DG_{\la}(\tilde c_j(\la))$ is pluriharmonic on 
$\Omega \setminus A$.
As it is continuous on $\Omega$, it is actually pluriharmonic on $\Omega$ and 
therefore coincides 
with
$\log|h(\la)|$ for some non-vanishing holomorphic function $h$.

\noindent \label{omega B}It remains to show that the coefficients of $H_n$ are uniformly bounded, for this choice of $h$.
Let us  consider an arbitrarily small ball $B$ contained in 
$\Omega \setminus 
A$. We will show that,
for all $\la \in  B$, one has $H_n(\la,z)=e^{-id^n\theta_B}\Pi (z 
\wedge A_j(\la))$
where $\theta_B \in {\bf R}$ and $A_j(\la)\in \{G_{\la}=0\}.$ The 
conclusion will follow since 
$\cup_{\la \in \Omega}\{G_{\la}=0\} \subset \subset {\bf C}^2$. As each term 
in the sum 
$\sum_{j=1}^DG_{\la}(\tilde c_j(\la))$ is {\em p.s.h.} on $B$, there are 
$D$ non-vanishing 
holomorphic
functions $h_j$ such that $G_{\la}(\tilde c_j(\la))=\log|h_j|$. Thus 
$\log|h|=\log\Pi|h_j|$ and 
$h=e^{i\theta_B}\Pi h_j$ for some $\theta_B\in \bf R$. Then, for 
any $\la \in 
B$, we get
from \eqref{zd}:
$$H_n(\la,z)=e^{-id^n\theta_B}\Pi\left(h_j(\la)^{-d^n}z\wedge F_{\la}^n
(\tilde 
c_j(\la)\right)=
e^{-id^n\theta_B}\Pi\left(z\wedge F_{\la}^n \left( \frac{\tilde c_j
(\la)}{h_j(\la)} \right)  
\right).$$
It finally suffices to set $A_j(\la):=\frac{\tilde c_j(\la)}{h_j(\la)} $ 
since, as desired, we have
$G_{\la}(A_j(\la))=G_{\la}(\tilde c_j(\la))-\log|h_j(\la)|=0$.

\bigskip

$\underline{Supp(dd^cL)\subset {\cal R}^c}$: this is given by Theorem \ref{a}.

$\underline{{\cal R}^c \subset \overline{\cal Z}_1(X,E
)}$: we shall use the 
following Lemma (\cite{BW}, 
Lemma VII.5).

\begin{lem}\label{ze} Let $z_0 \in \pp$ be a repulsive fixed point of $f^{n_0}_
{\la_0}$ ($n_0$ being 
the period of the associated
cycle) and $B$ be a ball centered at $\la_0$ in $X$. Let $z(\la)$ be a holomorphic map defined on 
some neighbourhood of
 $\la_0$ in $X$ such that $z(\la_0)=z_0$ and, for every $\la$,   $z(\la)$ is a repulsive fixed 
point of $f_{\la}^{n_0}$ 
(the points $z(\la)$
are given by the implicit function theorem). Then: 
either

\begin{itemize}\item[] i) $z(\la)$ holomorphically extends to $B$ and $z(\la)$ 
is a repulsive fixed 
point of  
$f_{\la}^{n_0}$ which belongs to a cycle of period $n_0$, for every $\la \in B$,\end{itemize}

or

\begin{itemize}\item[] ii) $z(\la)$ holomorphically extends
   to a neighbourhood of   some path 
$\gamma$ joining 
$\la_0$ to $\la_1$
in $B$ and $z(\la_1)$ is an attracting fixed point of $f_{\la_1}^{n_0}$. In 
particular, there are 
infinitely many values of
$\la'$ such that $z(\la')$ is a neutral fixed point of $f_{\la'}^{n_0}$ and 
the   set of   corresponding multiplier  s   
contains an open subset of $S^1$. Again every $z(\la)$ belongs to a cycle of 
period 
$n_0$.
 \end{itemize}
 \end{lem}

\bigskip
\noindent If $\la_0\in {\cal R}^c$ then, using the above lemma, we may find a 
non stationary 
sequence $\la_k \to \la_0$ such that $f_{\la_k}$
has a neutral cycle of period $n_k$ and a multiplier $e^{2i\pi 
\theta_k}$ with $\theta_k \in E$ . Since by Fatou 
theorem $f_{\la_0}$
has at most $6d-6$ non-repulsive cycles, all but a finite number of $Per
(X,n_k,e^{2i\pi \theta_k})$
differ from $X$, this shows that $\la_0 \in \overline{\cal Z}_1(X,E
)$.

\bigskip

$\underline{\overline{\cal Z}_1(X,E
) \subset {\cal S}^c}$: we proceed by 
contradiction. Let   $\la_0 \in 
\overline{\cal Z}_1(X,E
)$   and $B$ an open ball in $X$ such 
that $\la_0 \in B \subset {\cal S}$. Let $n_0 \in {\bf N}^*$ and $\theta_0 \in 
E$ such 
that 
$B \cap Per_1(X,n_0,e^{2i\pi \theta_0}) \neq \emptyset$. On a small ball $B' 
\subset B$ centered at 
some point 
$\la_1 \in Per_1(X,n_0,e^{2i\pi \theta_0})$ there exists a holomorphic map $z
(\la)$ such that 
$f_{\la}^{n_0}(z({\la}))=z({\la})$.
Moreover, as the multiplier of $f_{\la}^{n_0}$ at $z({\la})$ is not constant 
near $\la_1$ (otherwise 
  $Per(X,n_0,e^{2i\pi \theta_0})$   would coincide with $X$), we may find 
$\la_2,\la_3\in B'$ such 
that
$z({\la_2})$ (resp. $z({\la_3})$) is attractive (resp.repulsive) for $f_{\la_2}
^{n_0}$ (resp. 
$f_{\la_3}^{n_0}$). As the basin of 
$f_{\la_2}^{n_0}$ at $z({\la_2})$ contains a critical point, there exists $0 
\leq i_0 \leq n_0$ such that 
the sequence
$d[z({\la}),f_{\la}^{kn_0-i_0}(C_{f_{\la}})]$ is converging to $0$
 around $\la_2$. Then, since $B' \subset B \subset {\cal S}$, $d[z({\la}),f_
{\la}^{kn_0-i_0}(C_{f_{\la}})]$
actually converges to $0$ in $B'$ which is impossible because $z(\la_3)$ is 
repulsive. \qed

\begin{rem}\label{zg} {\em Two fundamental facts in Ma\~n\'e-Sad-Sullivan theory are 
the density of ${\cal 
R}$ in $X$
and the emerging concept of hyperbolic component: two elements lying in the 
same connected 
component of ${\cal R}$
are either both hyperbolic or both non-hyperbolic. This plays an important 
role in the approach of 
Fatou's conjecture 
on the density of hyperbolic rational maps. It turns out that these facts may 
be established by mean 
of elementary
 arguments similar to those used in the last steps of the proof of Theorem \ref
{zb}.}
\end{rem}

\bigskip

We now aim to generalize Theorem \ref{zb} to the case of powers $(dd^cL)^p$. 
To 
this purpose we 
have to discuss the intersection of
$p$ hypersurfaces $Per_1(X,n,e^{2i\pi \theta})$. For any $N_p:=(n_1,\dots,n_p) 
\in ({\bf N}^*)^p$ 
and 
$\Theta_p:=(\theta_1,\dots,\theta_p)\in E^p$ we define 
$$Per(X,N_p,e^{2i\pi \Theta_p}):=Per(X,n_1,e^{2i\pi \theta_1})\cap \dots \cap Per
(X,n_p,e^{2i\pi \theta_p}).$$
As previously,   $Per_p(X,N_p,e^{2i\pi \Theta_p})$ denotes the union of all the codimension $p$, irreducible
components   of
$Per(X,N_p,e^{2i\pi \Theta_p})$. We then set:
$${\cal Z}_p(X,E
):=\bigcup_{N_p\in ({\bf N}^*)^p,\Theta_p \in E^p}Per_p(X,N_p,e^{2i\pi \Theta_p}).$$

\bigskip
Our generalization may be stated as follows.

\begin{theo}\label{zf} Let $E$ be a dense subset of $]0,1[$. Let $\fL$ 
be a holomorphic family of rational maps of 
degree $d$ on $\pp$.
 Then for any $1 \leq p \leq dim_{\bf C}X$: 
$$Supp(dd^cL)^p \subset \overline{\cal Z}_p(X,E
).$$
\end{theo}

\proof We proceed by induction on $p$. Let us call $({\cal H})_p$ the 
following assertion:\\

{\em $({\cal H})_p$: For \underline {any} complex manifold $X$ of dimension $n 
\geq p$ and any 
holomorphic family
$\fL$ parametrized by $X$ we have $Supp(dd^cL)^p \subset \overline{\cal Z}_p
(X,E
).$ }\\ 

According to Theorem \ref{zb} $({\cal H})_1$ is true. Let us show that 
$({\cal H})_p$ implies $({\cal 
H})_{p+1}$.
To this end we shall combine the following fact with $({\cal H})_1$.\\

\noindent{\bf Fact:} {\em Assume that $({\cal H})_p$ is true. Let $U$ be an 
open set in ${\bf C}^n$ 
($n>p$) and
$\{f_{\la}\}_{\la \in U}$ be a holomorphic family. If $L(f_{\la})$ is 
pluriharmonic on every
$Per_p(U,N_p,e^{2i\pi \Theta_p})$ then $(dd^cL)^{p+1}\equiv 0$ on $U$.}

\bigskip

This fact, which is actually      the heart of our proof, will be established 
later. It is useful to remark 
that a continuous function 
on an analytic set $Y$ is {\em p.s.h.} if and only if
it is {\em p.s.h.}  on the set of regular points of $Y$ (see \cite{De}, 
Theorem 1.7). We also recall that $L$ is continuous (see Corollary \ref{n}).

Let us consider a holomorphic family $\fL$ and $\la_0\in Supp(dd^cL)^{p+1}$. 
Pick an arbitrarily 
small open set $U$
such that $\la_0 \in U$. We have to show that ${\cal Z}_{p+1}(X,E
)\cap U \neq 
\emptyset$. We may 
identify $U$ with an open set of 
${\bf C}^n$. According to the above fact, there exist $N_p \in ({\bf N}^*)^p$ 
and $\Theta_p \in E^p$ such
that $L$ is not pluriharmonic on $Per_p(U,N_p, e^{2i\pi \Theta_p})$. This 
implies the existence of 
some
regular curve $\Gamma$ contained in $Per_p(U,N_p,e^{2i\pi \Theta_p})$ such 
that 
$dd^c(L_{|\Gamma})$ does not
vanish. Thus, Theorem \ref{zb} applied to the family $\{f_{\la}\}_{\la \in 
\Gamma}$ guarantees the 
existence of some 
$Per_1(\Gamma,n_{p+1}, e^{2i\pi \theta_{p+1}})$,   $\theta_{p+1}\in E$.  \label{NB 10}Let $N_{p+1}:=(N_p,n_
{p+1})$ and
$\Theta_{p+1}:=(\Theta_p,\theta_{p+1})$; since $Per_{p+1}(U,N_{p+1},e^{2i\pi \Theta _{p+1}}) \subset {\cal Z}_{p+1}(X,E
)\cap U$,  it is enough to observe that $Per_{p+1}(U,N_{p+1},e^{2i\pi \Theta _{p+1}}) \neq \emptyset$. \qed
 
 . 

\bigskip

\noindent \underline{Proof of the fact}: By an elementary slicing argument, 
the positive current
$(dd^cL)^{p+1}$   vanishes   identically  on $U$ as soon as the
positive measures obtained by restriction on the $(p+1)$-dimensional affine 
subspaces  
  vanish.  
Let $S$ be the intersection of $U$ with such an affine subspace of ${\bf C}^n$ 
and
set $L_0:=L_{|S}$, $\mu_0:=(dd^cL_0)^{p+1}$. We have to show that
for every euclidean, $(p+1)$-dimensional, open ball $B \subset \subset S$,
the measure $\mu_0$ vanishes on $\frac{1}{2}B$. 

To this end, we introduce the solution $\tilde L_0$ of the Dirichlet-Monge-Amp\`ere problem
with   datum   $L_0$ on $bB$. The function $\tilde L_0$ is continuous on $\overline 
B$, coincides
with $L_0$ on $bB$ and is {\em p.s.h.} maximal on $B$ (see \cite{BT}). By 
maximality,
$\tilde L_0 \geq L_0$ on $\overline B$. We also consider the set $\Sigma_
{\varepsilon} \subset 
\frac{1}{2}B$
where $L_0$ and $\tilde L_0$ are $\varepsilon$-close:
$$\Sigma_{\varepsilon}:=\{0 \leq \tilde L_0 -L_0 \leq \varepsilon\} \cap \frac
{1}{2}B.$$
A theorem of Briend-Duval (see \cite{BD} or \cite{Sb} Theorem A.10.2) states 
that
$$\mu_0(\Sigma_{\varepsilon}) \leq C\varepsilon$$
where the constant $C$ depends only on $L_0$ and $B$. 
It thus suffices to show that $Supp (\mu_0) \cap \frac{1}{2}B \subset \Sigma_
{\varepsilon}$ for any $\varepsilon >0$.\\
 
The set ${\cal Z}_p(S,E
)$ is an union of complex curves in $S$. Let $A$ be one 
of these curves,
 then $A$ is a component of  
$Per_p(S,N_p,e^{2i\pi \Theta _p})$ and is therefore contained  in $S \cap Per_p
(U,N_p, e^{2i\pi \Theta _p})$; since $S$ is an affine subspace this is easy to check on the regular part of $A$. Thus the 
function $L=L_0$ is,
by assumption, harmonic on $A\cap B$. The function $\tilde L_0-L_0$ is thus 
subharmonic on 
$A\cap B$
and, by the maximum principle,   vanishes   identically. Therefore $\tilde L_0-L_0
$ vanishes
on ${\cal Z}_p(S,E
)\cap B$. 

Of course $Supp(\mu_0) \subset Supp(dd^cL_0)^p$ and,  as  $({\cal H}_p)$ is 
supposed be true,  
\linebreak $Supp(dd^cL_0)^p \subset \overline {\cal Z}_p(S,E
)$. Thus  $\tilde 
L_0-L_0$ vanishes
on $Supp(\mu_0) \cap \frac{1}{2}B$, which is, therefore, contained in $\Sigma_
{\varepsilon}$.  \qed

\bigskip

\section{  The bifurcation measure  }
In this section  we define
 the bifurcation measure $\mu$ on the moduli space ${\cal M}_d$ of rational maps $\pp \to \pp$ 
and we establish some basic results about it. Although the section is mainly devoted to the one-dimensional case, 
the fact that the bifurcation currents
$(dd^cL)^p$ have finite mass on $\Hd$ will be established in any dimension, that is in $\Hd(\PP)$.

\bigskip
The group $PSL(2,{\bf C})$ of M\"obius transformations acts on the space $\Hd(\pp)$
by conjugation. Two conjugated rational functions $f_1,f_2\in {\cal H}_d(\pp)$ 
 have the same dynamics, therefore in order to study the stability of holomorphic
families  of rational functions, one
 considers, instead of ${\cal H}_d(\pp)$, the moduli space
${\cal M}_d:={\cal H}_d({\bf P}^1)/PSL(2,{\bf C})$. 

\begin{rem}\label{ai} {\em The moduli space ${\cal M}_d$ is a normal, quasi-projective variety (see \cite{S}, Remark p.43); 
the proof requires some effort because 
 $PSL(2,{\bf C})$ is not compact and its action on $\Hd(\pp)$ is not free (indeed
  there is some $f\in {\cal H}_d({\bf P}^1)$  whose   isotropy group
$Aut(f):=\{\varphi \in PSL(2,{\bf C}); \varphi^{-1}\circ f \circ \varphi =f \}$
 is not trivial). 
Here we recall some useful facts about ${\cal M}_d$: 
}\end{rem}

\begin{itemize}
\item[(i)] {\it the canonical projection $\Pi:{\cal H}_d({\bf P}^1) \to {\cal M}_d$ is open;}
\item[(ii)]{\it for all $f\in {\cal H}_d({\bf P}^1)$, the isotropy group $Aut(f)$ 
is finite and locally 
there is a complex submanifold $V$, invariant by
the action of $Aut(f)$,
 transverse at $f$ to the orbit of $f$,  such that $\Pi (V)$ is open in ${\cal M}_d$ and the canonical projection $\Pi$ induces a biholomorphism
$V/Aut(f) \to \Pi(V)$;}
\item[(iii)] {\it the set of all $f \in {\cal H}_d({\bf P}^1)$ such that
 $Aut(f) \neq \{id_{{\bf P}^1}\}$ is an analytic subset $Z$ of
${\cal H}_d({\bf P}^1)$ and} Sing$({\cal M}_d) \subset \Pi(Z)$.\end{itemize}

\medskip

\noindent(In general Sing$({\cal M}_d) \neq \Pi(Z)$, e.g.  ${\cal M}_2={\bf C}^2$ is smooth, 
and $\Pi(Z)$ is a cubic curve of ${\bf C}^2 = {\cal M}_2$ (see \cite{M1}, Corollary 5.3).
But, for $d > 2$, ${\cal M}_d$ has singular points). It follows that if $f\notin Z$,
then  $PSL(2;{\bf C})\times V \simeq \Pi^{-1}(\Pi(V))$, therefore 

\begin{itemize}
\item[(iv)]
${\cal H}_d({\bf P}^1) \setminus Z \to {\cal M}_d \setminus \Pi(Z)$
{\it is a principal bundle.}\end{itemize}

\medskip
\noindent Since $\dim_{\bf C} PSL(2,{\bf C})=3$ and all isotropy groups $Aut(f)$ are
finite,  hence each orbit is a complex $3$-fold and
$\dim_{\bf C}{\cal M}_d=\dim_{\bf C}{\cal H}_d({\bf P}^1)-3=2(d-1)$. The 
 Lyapunov exponent $L(f)$,  $f \in \Hd(\pp)$, which is invariant
  under   the action of $PSL(2,{\bf C})$, is constant on the orbits and, if
$p> 2(d-1)$,
 the current $(dd^cL)^p$ vanishes identically 
on $\Hd(\pp)$. Therefore in order to define a measure   by means of Monge-Amp\`ere operator on $L$,
it is necessary to consider the function
   $\hat L:{\cal M}_d \to {\bf R}$ induced from $\Hd(\pp) \stackrel{L}{\longrightarrow} \bf R$.

\begin{prop}\label{aj} The function $\hat L$ is continuous, bounded from below and {\em p.s.h.} 
on ${\cal M}_d$.\end{prop}

\proof  By \cite{R}, $L$ is bounded from below. Using Corollary \ref{n} it is enough to   notice  
 that, by means  of Remark \ref{ai} (iv), $\hat L$
is {\em p.s.h.} on ${\cal M}_d \setminus \Pi(Z)$ and thus (see \cite{De}, Theorem 1.7) on the whole ${\cal M}_d$. \qed

\bigskip

Now the  currents $(dd^c\hat L)^p$, $1 \leq p \leq 2(d-1)$, are well defined on ${\cal M}_d$. 
In particular,    the measure $\mu:=(dd^c\hat L)^{2(d-1)}$ will be called {\em bifurcation measure}.  

\medskip

\begin{prop}\label{al} The   bifurcation measure   $\mu$ does not vanish identically, in particular any
non-flexible Latt\`es map lies in the support of $\mu$.
\end{prop}

\proof  For $d\geq 2$ fixed, let 
$f_0\in {\cal H}_d({\bf P}^1)$  be a non-flexible Latt\`es map 
(e.g. $f_0$ is the map associated to an imaginary quadratic number field, see \cite{M}
 Lemma 5.4), then all  Latt\`es maps which belong to a small neighbourhood
of $f_0$ in $\Hd(\pp)$ are conjugated to $f_0$.
Let $V$ be a complex submanifold  in a neighbourhood of $f_0$ as in Remark \ref{ai} (ii); 
since the function $L(f)$ takes its minimum value $\log \sqrt{d}$ exactly when $f$
is a Latt\` es map (see \cite{Le}, \cite{Zd}), hence $f_0$ is a point of strict minimum for $L_{|V}$.
As  $dim_{\bf C}V=2(d-1)$ we shall see that $f _0\in Supp(dd^cL_{|V})^{2(d-1)}$,
i.e. $\Pi(f_0)\in Supp(\mu)$.  For every  small, euclidean, open ball $B \subset V$ centered at $f_0$
there is a suitable constant $c$ such that $L(f_0) < c < L(f)$, for every $f \in B$; so the function $L_{|V}-c$ does not take
its minimum on $\overline B$
at the boundary, therefore (see  Theorem A in \cite{BT}) it is not maximal, that is $(dd^cL_{|V})^{2(d-1)}$
does not vanish identically on $B$. \qed

\bigskip

In order to see that  $\mu$ has finite mass (see Proposition \ref{ar}) we shall show that $L$
extends from $\Hd$ to the whole projective space across the hypersurface $\Sigma_d$ and that the powers of $dd^cL$ have finite mass
on  ${\cal H}_d$.
Since these results hold for holomorphic maps ${\bf P}^k \to {\bf P}^k$, $k \geq 1$, we believe that it is useful
to present them in this more general case.

First of all let us recall that $\Hd({\bf C}^{k+1})={\bf C}^N \setminus \tilde \Sigma _d$, see Subsection \ref{sssp} and 
that from  \eqref{o1} 
and Proposition \ref{z} it follows:
\begin{equation}\label{am}L(F)=H(F)-\frac{k+1}{d^k}\log|Res(F)|+cst.\end{equation}
for every polynomial map  $F\in \Hd({\bf C}^{k+1}) \subset {\bf C}^{N+1}$.

\begin{prop}\label{an}
The function $H$ extends from $\Hd({\bf C}^{k+1})$ to a {\em p.s.h.} function on the whole ${\bf C}^{N+1}$ 
and  the function $L(f)$ extends from ${\cal H}_d({\bf P}^k)$ to a function $L^1_{loc}({\bf P}^N)$. Moreover
there is a $(1,1)$-current $R$, 
positive and closed on ${\bf P}^N$ such that $\pi^*R=dd^cH$
 and
\begin{equation}\label{ao}dd^cL=R-\frac{k+1}{d^k}[\Sigma_d].\end{equation}\end{prop}

\proof  Let $ {\cal K} \subset {\bf C}^{N+1}$ be compact, let us check that $H$ is bounded from above
on ${\cal K} \cap {\cal H}_d({\bf C}^k)$. By definition  (see \eqref{o2}):
$$H(F)=\int_{\PP}g_F\sum_{j=0}^{k-1}T_f^j\wedge\omega^{k-j-1}\wedge [C_f]+\int_{{\bf P}^k}\log||J_F||_0\omega^k.$$
Then, as $g_F$ is locally bounded from above (see Remark \ref{1aa}), one concludes taking into account the following formulas:
  $\int_{\PP}\sum_{j=0}^{k-1}T_f^j\wedge\omega^{k-j-1}\wedge [C_f]=k \ deg(C_f)$ and  $\int_{{\bf P}^k}\log||J_F||_0\omega^k=$ 
 $=\int_{{\bf C}^{k+1}}\log|det(F'(z))|m$ (Lemma \ref{6}).
From this, since $H$ is {\em p.s.h.} on $\Hd({\bf C}^{k+1})$ (see the Claim  in the proof of Theorem \ref{s}), it follows that
$H$ extends to a {\em p.s.h.} function on the whole  ${\bf C}^{N+1}$.  Thus the right hand side of
\eqref{am} belongs to $L^1_{loc}({\bf C}^{N+1})$ and extends $L(F)$ as a $0$-homogeneous $L^1_{loc}$ function
on ${\bf C}^{N+1}$. Thus $L(f)$ is well defined on the whole ${\bf P}^N$. Choosing 
a holomorphic section $U \stackrel{\sigma}{\longrightarrow} {\bf C}^{N+1}\setminus\{0\}$ on an open subset
 $U$ of ${\bf P}^N$, we get
 $\forall f \in U$, 
\begin{equation}\label{ap}L(f)=H(\sigma (f))-\frac{k+1}{d^k}\log|Res(\sigma (f))|
+cst.\end{equation}
As $dd^c(H \circ \sigma)$ does not depend on $\sigma$, it defines a positive, closed current $R$ on ${\bf P}^N$
such that  $\pi^*R=dd^cH$.  Then   
\eqref{ao}  follows from \eqref{ap}. \qed  

\bigskip
The sum $L(f)$ of the Lyapunov exponents is bounded from below 
(see  Theorem 1 in \cite{BD}), thus, as in the one-dimensional case, the powers
of $dd^c(L_{|{\cal H}_d({\bf P}^N)})$ are well defined; but
 to show that these currents have finite mass
  requires   some work.

\begin{prop}\label{aq} For $1 \leq p \leq N$, $$\int_{\Hd(\PP)}(dd^cL)^p\wedge \omega^{N-p} < \infty$$ and 
the trivial extension $S_{(p)}$ of $(dd^c(L_{|{\cal H}_d}))^p$ to the whole ${\bf P}^N$ is well defined.
\end{prop}

\medskip

\noindent{\bf Remark.} We recall that, by definition of {\em trivial extension}, $S_{(p)}$  is the positive, closed, $(p,p)$-current on ${\bf P}^N$
characterized by
\begin{itemize}\item[(i)]
$S_{(p)}=(dd^cL)^p$ on ${\cal H}_d({\bf P}^k)$ 

\item[(ii)]  $\chi_{\Sigma_d}.S_{(p)}=0$, where $\chi_{\Sigma_d}$ is the characteristic function.
\end{itemize}
\proof From \eqref{ao} it follows $dd^cL \leq R$, thus $dd^cL$ has finite mass on $\Hd(\PP)$ and
its trivial extension to ${\bf P}^N$ is the positive, closed current
$S_{(1)}:=(1-\chi_{\Sigma_d})R$.

Now we shall argue by induction, assuming that the trivial extension $S_{(p)}$ of
$(dd^c(L_{|{\cal H}_d}))^p$ to the whole ${\bf P}^N$ is well defined 
(and, of course, positive and closed).
There is a smooth, closed $(1,1)$-form $\alpha$ on ${\bf P}^N$ such that $S_{(1)}-\alpha =dd^cu$ and, 
by means of the regularization theorem of 
Demailly (\cite{De1}), there are 
a sequence $ \{u_n \}$ of smooth functions decreasing to $u$ and a sequence  $\{\la_n\}$ of continuous
functions decreasing to
$\nu (S_{(1)},.)$, such that $S_n:=\alpha + dd^cu_n \to S_{(1)}$ and $S_n +\la_n \omega \geq 0$.

We can estimate the mass of $S_{(p)}\wedge (S_n +\la_n \omega)$ (which is a positive current)
as follows:
$$||S_{(p)} \wedge (S_n+\la_n\omega)||=\int_{{\bf P}^N}S_{(p)}\wedge \alpha \wedge \omega^{N-p-1}+$$
$$+\int_{{\bf P}^N}
S_{(p)} \wedge dd^cu_n \wedge \omega^{N-p-1}+\int_{{\bf P}^N}\la_nS_{(p)}\wedge \omega^{N-p}$$
Let us look at the right hand side: the first term is constant, the second vanishes. Since $L$ (and therefore $H$)
is bounded from below on 
${\bf P}^N \setminus \Sigma_d=\Hd(\PP)$,   the   Lelong numbers of  $S_{(1)}$ vanish outside $\Sigma_d$. Thus
$\la _n$ decreases to $\nu(S_{(1)},.)\chi_{\Sigma_d}$, but $\chi_{\Sigma_d}S_{(p)}=0$, therefore $\la_nS_{(p)}\to 0$.

This means that  $S_{(p)} \wedge (S_n +\la_n \omega)$ has bounded  mass,  thus  
we can assume that it converges to a positive current $Q$. But $S_{(p)}\wedge \la_n \omega \to 0$, thus $Q$ is closed.
 Therefore we can set
$S_{(p+1)}:=(1-\chi_{\Sigma_d})Q$. \qed

\bigskip

Coming back to one-dimensional case, we can establish:

\begin{prop}\label{ar} The bifurcation measure $\mu$ has finite mass on ${\cal M}_d$.\end{prop}

\proof Since $\hat L$ is bounded from below, the measure $\mu$ does not charge analytic subsets 
thus from Remark \ref{ai} (iv) and the previous Proposition we get:
$$\int_{{\cal M}_d}\mu=\int_{{\cal M}_d \setminus \Pi(Z)}(dd^c\hat L)^{2(d-1)}=$$
$$=\int_{{\cal H}_d({\bf P}^1)\setminus Z}
(dd^c L)^{2(d-1)}\wedge \omega^3=\int_{{\cal H}_d({\bf P}^1)}(dd^c L)^{2(d-1)}\wedge \omega^3 < \infty.$$
\qed

\bigskip

\begin{rem}\label{as} {\em Since the hypersurfaces $Per(\Hd (\pp),n,e^{2i\pi \theta})$ are invariant under the action of $PSL(2,{\bf C})$,
there are no difficulties in order to get, from Theorems \ref{zb} and \ref{zf}, the corresponding statement for $\hat L$.} Actually
the following   claim   holds: $Supp(dd^c \hat L)=\overline{\cal Z}_1({\cal M}_d)$ and for $1 < p \leq 2(d-1)$,
$  Supp(dd^c \hat L)^p \subset \overline{\cal Z}_p({\cal M}_d)$. In particular 
$Supp(\mu) \subset \overline{\cal Z}_{2(d-1)}({\cal M}_d)$
\end{rem}
\bigskip

We say that a point $x\in {\cal M}_d$ is {\em chaotic} if the Julia set of any $f\in \Pi^{-1}(x)$ is $\pp$.  

\begin{prop}\label{at} In any neighbourhood of a point of $Supp(\mu)$ there are uncountably many chaotic points.
\end{prop}

\noindent \underline{Proof}\footnote{We would like to thank T.C. Dinh who told us the possibility
to use Shishikura's theorem here.}. Consider  $\{\theta \in {\bf R}; \limsup_{n\to \infty}\frac{\log \log(1/|\theta^n-1|)}{n}<\log d\}$ and use this
open dense and uncountable subset of $\bf R$ to define $E$ and ${\cal Z}_p(\Hd (\pp),E)$ (see \cite{M2}). If $z_0$ is a periodic points of
$f\in \Hd (\pp)$ with
multiplier $e^{2i\pi \theta}$, $\theta \in E$, then $z_0$ is a Cremer point; therefore any $f \in {\cal Z}_{2(d-1)}(\Hd (\pp),E)$
has $2(d-1)$ Cremer points and (see \cite{Shi}, Corollary 2) is chaotic. \qed

\section{Examples and applications}

\subsection{Geodesics on the space of  K\" ahler metrics}\label{S7a}
Let  $M$ be a compact K\" ahler manifold of dimension $k$ with a fixed  K\" ahler metric $\omega$. 
In order to discuss extremal (e.g. Einstein and of constant scalar curvature) K\"ahler metrics 
it is useful
to consider the space $H_{\omega}$ of K\" ahler metrics with the same K\" ahler class of $\omega$ (see e.g. \cite{Chen}). 
It can be   thought   also as the space
of K\" ahler potentials, that is
$$H_{\omega}:=\{\phi \in {\cal C}^{\infty}(M); \omega +i \partial \overline \partial \phi > 0 \}/\sim $$ 
where $\phi_1 \sim \phi_2$ if and only if $\phi_1-\phi_2$ is a constant. Endowing $H_{\omega}$ with a suitable metric it turns out 
that $H_{\omega}$, 
as Riemannian manifold, is an infinite dimensional
symmetric space and there is a (unique) Levi-Civita connection whose curvature is covariant constant (see \cite{Sem} and \cite{Chen}).

For such a connection  the equation of geodesic is
\begin{equation}\label{ab1}(i\partial \overline \partial \phi + \omega)^{k+1}=0\end{equation}
This means that $\phi$ is a smooth real function  defined on $[0,1] \times M$ (in this case one understands the $\overline \partial$ operator 
as the one on the cylinder $[0,1] \times S^1$ with its natural complex structure) or, for complex geodesics, it is more generally defined on
$X \times M $ where $X$ is a
Riemann surface. We point out that very few explicit examples of these geodesics are known, thus the following remark
may have some interest.

Let $M=\PP$ and $\omega$ be the Fubini-Study metric, Corollary \ref{v}  says that any holomorphic family $\fL$ of endomorphisms of
$\PP$ defines a 
``geodesic'' $\phi:=g_{F_{\la}}$. Of course the behaviour of $dd^cg_{F_{\la}}+\omega$ is very far from
the desired regularity, but  there is at least one case
in which holomorphic dynamics may give interesting examples: let $M={\bf P}^1$ and $\fL$ be a family of flexible 
Latt\`es maps (see \cite{M}, Ch. 8.3), then the functions $g_{F_{\la}}:{\bf P}^1 \to \bf R$ are smooth outside a
finite set.

\subsection{Attractor  s   in ${\bf P}^2$}\label{S7b}

\begin{defi}\label{af} Let $\{f_{\la}\}_{\la \in X}$ be an one parameter holomorphic family
 (i.e. $X$ is an open subset of ${\bf C}$) of
endomorphisms of ${\bf P}^k$ and let $Y$ be a complex subspace of $X\times {\bf P}^k$ of
pure dimension $q$.
We shall say that the Green function $G_{\la}$ is {\em maximal} on $Y$ if and only if,
for every holomorphic section
$\PP \supset U \stackrel{\sigma}{\longrightarrow} {\bf C}^{k+1}\setminus\{0\}$, 
$$\big(dd^c(G_{\la} \circ \sigma)\big)^q=0\  \hbox{on} \ Y \cap (X \times U).$$\end{defi}

\noindent Although the Green function depends on the choice of the lift of $f_{\la}$, 
the definition is well posed since 
$dd^c(G_{\la} \circ \sigma)$ does not
depend on the particular family of lifts $\{F_{\la}\}_{\la \in W}$
chosen in an open subset $W$ of $X$.

With this definition we can give the following formulation of Corollary \ref{v}

\begin{prop}\label{ag} Let $\{f_{\la}\}_{\la \in X}$ be an one 
parameter holomorphic family  of
endomorphisms of ${\bf P}^k$. Then $G_{\la}$ is maximal on $X \times {\bf P}^k$.
Moreover $G_{\la}$ is maximal on $C_X$ if and only if $L(f_{\la})$ is harmonic.
\end{prop}
\bigskip
Now we shall apply this Proposition to a particular case. 
 For $\varepsilon \in {\bf C}$, consider the rational map ${\bf P}^2 \to {\bf P}^2$ defined by
$$f_{\varepsilon}=[P(z,w):Q(z,w):t^d+\varepsilon R(z,w)]$$
where $P,Q,R$ are homogeneous polynomials of degree $d \geq 2$ such that $(P,Q)$
is non degenerate and the induced rational function $$f=[P(z,w):Q(z,w)]$$
is strictly critically finite. It is useful to consider the line 
${\cal R}_{\infty}:=\{t=0\}$ as the line
at infinity of the complex plane ${\bf C}^2 \simeq \{[z:w:1] \in {\bf P}^2 \}$. 
If $a\in {\cal R}_{\infty}$, we shall denote by ${\cal R}_a$ the line passing through
the origin $[0:0:1]$ and $a$.
The map
$f_{\varepsilon}$ preserves lines through the origin and moves them in a chaotic way,
since $f$ is chaotic, indeed identifying ${\cal R}_{\infty}$ with ${\bf P}^1$ we get 
$f_{\varepsilon}({\cal R}_a)\cap {\cal R}_{\infty}=\{f(a)\}. $

For $|\varepsilon|<<1$ the only 
Fatou component is the superattractive basin of the origin 
(see \cite{FS-Ex} Lemma 2.1).  Moreover 
for $|\varepsilon|<<1$ the map $f_{\varepsilon}$ has an attractor $A$ contained in a
neighbourhood of the line at infinity, which intersect any line passing through the origin (ibidem,
Lemma 2.2).

Our aim is to show that, for $|\varepsilon|<<1$, the family $\{f_{\varepsilon}\}$ is stable
in the following sense:

\begin{prop}\label{ah}
 The function $\varepsilon \to L(f_{\varepsilon})$ is
harmonic in a neighbourhood of $0 \in \bf C$.
\end{prop}

\proof Let $X:=\bf C $
and $F_{\varepsilon}=(P,Q,t^d+\varepsilon R)$. A simple inspection on $F'_{\varepsilon}$
shows that 
$C_{f_{\varepsilon}}={\cal R}_{\infty} \cup \big( \cup_{c \in C_f}{\cal R}_c \big)$,
thus the critical set does not depend on $\varepsilon$ and
$$C_X=\big(X \times {\cal R}_{\infty}\big) \cup \big( \cup_{c \in C_f}X \times{\cal R}_c \big).$$

Let $c\in C_f$; then by hypothesis, there exist $j,k \in {\bf N}^*$ such that, 
$a:=f^j(c)$ and $f^k(a)=a$. This means that ${\cal R}_a$ is fixed by $f_{\varepsilon}^k$.
Thus we can consider the family $f_{\varepsilon |{\cal R}_a}^k$ as a family
of endomorphisms of ${\cal R}_a$. From  Proposition \ref{ag}  it follows that, for this
family, the Green function is maximal on $X \times {\cal R}_a$. But all the powers 
of $F_{\varepsilon}$ have the same  Green function $G_{F_{\varepsilon}}$. That is,
for every section $U \stackrel{\sigma}{\longrightarrow} {\bf C}^{3}\setminus\{0\}$,
the function $G_{F_{\varepsilon}}\circ\sigma $ is maximal on
$X \times  (U\cap {\cal R}_a)$. Since $f_{\varepsilon}^j({\cal R}_c)={\cal R}_a$
and $G_{F_{\varepsilon}}\circ\sigma \circ f_{\varepsilon}^j=
G_{F_{\varepsilon}}(h. (F_{\varepsilon}^j \circ \sigma'))=d^jG_{F_{\varepsilon}}\circ \sigma'+
\log|h|$, for a suitable section $\sigma'$ and a holomorphic,
never vanishing, function $h$,  we get that $G_{F_{\varepsilon}}$ is maximal on
$X \times {\cal R}_c$.

Now we shall show that choosing $V=\{|\varepsilon|<<1\}$, the function 
$G_{F_{\varepsilon}}$ is maximal on $V \times {\cal R}_{\infty}$. Let 
$u \in {\cal R}_{\infty} \simeq {\bf P}^1$ be a periodic point for $f$
(that is $f_{\varepsilon}^j({\cal R}_u)={\cal R}_u$ for some $j$). From the proof
of Lemma 2.2 in \cite{FS-Ex} it follows that there is an open neighbourhood
$W$ of ${\cal R}_{\infty}$ in ${\bf P}^2$ such that, if $|\varepsilon|<<1$,
then $f_{\varepsilon}(W) \subset \subset W$; thus the family 
$\{f_{\varepsilon}^n(u)\}_{n \in {\bf N}}$ is a normal in $V$ (as functions of
$\varepsilon$). Therefore $\varepsilon \to G_{F_{\varepsilon}}(\sigma(u))$ 
is harmonic on $V$. That's enough since these points $u$ are dense in ${\cal R}_{\infty}$.
\qed

\subsection{The   bifurcation measure    on ${\cal M}_2$}
As we have already recalled the moduli space ${\cal M}_2$ of the rational functions of degree $2$ can be
identified biholomorphically and in a canonical way  with ${\bf C}^2$ (see \cite{M1} Remark 3.3). Such an identification involves 
the affine structure since, for
every $\eta \in {\bf C}$, $Per({\cal M}_2,1,\eta)$ is a straight line of ${\bf C}^2$. In particular, the Mandelbrot family
$\{z^2+c; c\in {\bf C}\}$ coincides with the straightline $Per({\cal M}_2,1,0)$ of rational functions with a superattractive fixed
point. 

\begin{prop}\label{A} The Mandelbrot family is disjoint from the support of the bifurcation measure
$\mu$ on ${\cal M}_2$.
\end{prop}
\proof 
Since every $f_0 \in Per({\cal M}_2,1,0)$ has
a superattractive fixed point, there is an open subset $V$,  $Per({\cal M}_2,1,0) \subset V \subset {\cal M}_2$,
 such that every $f \in V$ has an attracting fixed point.

For every $f \in {\cal M}_2$ the number of attracting or indifferent cycles 
 is $\leq 2$ (see \cite{Shi}, Corollary 1), therefore $V \cap {\cal Z}_2({\cal M}_2,E)= \emptyset$. 
Thus from Theorem \ref{zf} it follows that $V$ is disjoint from $Supp(dd^c\hat L)^2=Supp(\mu)$. \qed

\begin{rem}\label{B} {\em Let us point out that it is not necessary to use Shishikura's result. First we can assume that
if $f \in V$, then the attractive fixed point depends holomorphically on $f$. Consider 
a holomorphic disc $\{f_{\la}\}_{\la \in \Delta}$ contained in $V \cap Per_1({\cal M}_2,n,e^{2i\pi \theta})$.
If $\hat L$ is not harmonic on this disc, then from Theorem \ref{zb} there is a holomorphic function
$z(\la)$ on an open disc $\Delta' \subset \Delta$ such that $f_{\la}^n(z(\la))=z(\la)$ and two values
$\la_1,\la_2 \in \Delta'$ such that $z(\la_1)$ is repulsive and $z(\la_2)$ is attractive. Thus $f_{\la_2}$
has two attracting points, and  it is stable in ${\cal M}_2$. This contradicts the fact that
$f_{\la_2}\in Per_1({\cal M}_2,n,e^{2i\pi \theta}) \subset Supp(dd^c \hat L)$.
Therefore $\hat L$ is pluriharmonic on every $Per_1(V,n,e^{2i\pi \theta})$; from the Fact in the proof of
 Theorem \ref{zf} it follows that $(dd^c\hat L)^2=0$   on $V$.  }
\end{rem}

\section*{Appendix}

\subsection*{Proof of proposition \ref{t}}
Let $u$ be a locally bounded function on $X$; using \eqref{u} with $R=u$, 
Fubini theorem and a suitable partition of the unity,
it follows that $u_y=u_{|\pi^{-1}(y)}$, for a.e. $y \in Y$.
Since $T$ is positive and closed, then the slices of $T$ exist (see \cite{Siu} (10.3)); 
from 
\eqref{u},
since slicing commutes with the operators $d$, $\partial$ and $\overline \partial$  (see 
\cite{Siu} (10.4)), it follows that, for a.e. $y \in Y$,
$T_y$ is a current on $\pi^{-1}(y)$, positive and closed. Thus $u_yT_y$ is well defined. 
By definition $dd^cu \wedge T=dd^c(u \wedge T)$, thus we shall argue by recurrence and, in order to finish the proof,
 it is enough to show $u_yT_y=(uT)_y$. 

This obviously holds if $u \in {\cal C}^{\infty}(X)$. Fix $\phi$ and $\psi$ as in \eqref{u} and let $K$ be a compact set such that
$Supp(\psi)\subset K$.
For a.e. $y \in Y$
$$|\int_{\pi^{-1}(y)}u_yT_y\wedge \iota^*_y\psi| \leq
\sup_K|u|\  C_{\psi}\   ||T_y||_{K \cap \pi^{-1}(y)},$$
therefore 
$$|\int_Y\big(\int_{\pi^{-1}(y)}u_yT_y\wedge \iota^*_y\psi \big)\phi| \leq 
\sup_K|u|\  C_{\phi, \psi}\   ||T||_K$$
and the left hand side is well defined.

That's all, because this inequality shows that
 the operator $\Phi(u)=$ \linebreak
$=\int_Y \big( \int_{\pi^{-1}(y)} T_y \wedge \iota_y^*(u \psi) \big) \phi$ can be continuously extended 
from ${\cal C}^{\infty}(X)$ to  $L_{loc}^{\infty}(X)$. \qed

\subsection*{Proof of Proposition \ref{aa}} 
In order to compute the constant $C_{d,k}$ it is enough to consider a particular $F$; infact if we take 
$F(z_0, \dots, z_k):=(z_0^d,\dots ,z_k^d)$, then $Res(F)=1$ by Proposition \ref{x} and, from Proposition \ref{z}, it follows 
\begin{equation}\label{ac}C_{d,k}=B(F).\end{equation}
  It turns out that  \label{NB 8}
\begin{equation}\label{acc}G_F(z)=\log \max_{_{0 \le j \le k}}|z_j|.\end{equation}
Therefore $G_F$ and (see \eqref{o3}) $B(F)=C_{d,k}$ do not depend on $d$. We shall write $C_k:=C_{d,k}$.

Now a direct computation of $B(F)$ is possible using \eqref{o3} and \eqref{acc}, but the following proof is more elementary. 
   
\medskip

\noindent{\bf Claim 1:} $L(F)=(k+1)\log d.$

\noindent\underline{Proof of  Claim 1}. Using 1.4(i) it is enough to compute 
$det F'(z)=d^{k+1}z_0^{d-1}\dots z_k^{d-1}$ and remark that 
$Supp(\mu_{F})=\{z \in {\bf C}^{k+1}; |z_0|=\dots=|z_k|=1\}$. \qed

\medskip

\noindent{\bf Claim 2:}  $H(F)=(k+1)(d-1) \left( C_{k-1}-\frac{1}{2}(1+\frac{1}{2}+ \dots +\frac{1}{k})\right)+(k+1)\log d$

\noindent\underline{Proof of  Claim 2}. The critical set $C_f$ is the union of the projective hyperplanes
$H_s:=\{z_s=0\}$, $0 \leq s \leq k$,   more precisely $[C_f]=(d-1)\sum_{s=0}^k[H_s]$.
The restriction $\tilde F$ of  $F$ to the hyperplane $\{z_k=0\}$ is of the same form: 
  $\tilde F (z_0,\dots ,z_{k-1})=(z_0^d,\dots,z_{k-1}^d)$. Hence $B (\tilde F)=C_{k-1}$.
Now from \eqref{o3} it follows 
$$\int_{\PP}g_F\sum_{j=0}^{k-1}T_f^j\wedge\omega^{k-1-j}\wedge [C_f]=
(d-1)\sum_{s=0}^k\int_{H_s}g_F\sum_{j=0}^{k-1}T_f^j\wedge\omega^{k-1-j}=$$
\begin{equation}\label{ad}=(d-1)(k+1)B(\tilde F)=(d-1)(k+1)C_{k-1}.\end{equation}  
 From Lemma \ref{6} (2) 
$$\int_{{\bf P}^k}\log ||J_F||_0 \omega^k=\int_{{\bf C}^{k+1}}\log |detF'|m=
(k+1)\log d+(d-1)\sum_{j=0}^k\int_{{\bf C}^{k+1}}\log|z_j|dm,$$
and by means of an elementary computation
\begin{equation}\label{ae}\int_{{\bf P}^k}\log ||J_F(z)||_0 \ \omega^k=
(k+1)\log d-(k+1)(d-1)\frac{1}{2}(1+\frac{1}{2}+ \dots +\frac{1}{k}).\end{equation}
Putting \eqref{ad}-\eqref{ae} in \eqref{o2} we get the Claim. \qed

\medskip

Now putting \eqref{ac} and  the two claims in \eqref{o1}  we get
  $$C_k=C_{k-1}-\frac{1}{2}\sum_{j=1}^k\frac{1}{j}.$$  
Now to finish it is enough to   find   $C_1$; for $k=1$ the map $F$ is given by  $F(z_0,z_1)=(z_0^d,z_1^d)$ and in this case
\eqref{o1} gives $$L(F)=g_F(0:1)+g_F(1:0)+\int_{{\bf P}^1}\log ||J_F(z)||_0 \omega-2(d-1)B(F);$$
now by Claim 1, $L(F)=2\log d$, by \eqref{ae}, $\int_{{\bf P}^1}\log ||J_F(z)||_0 \omega=2\log d-(d-1)$,  and,   by \eqref{acc} 
  , $g_F(0:1)=g_F(1:0)=0$, thus:
$$-1/2=B(F)=C_1.$$ \qed

\bibliographystyle{amsalpha}

\footnotesize{Giovanni Bassanelli, Dipartimento di Matematica, Universit\`a di Parma, Parco Area delle Scienze, 

153/A - I- 43100 Parma, Italia.

{\em Email: giovanni.bassanelli@unipr.it}}
\bigskip

\footnotesize{Fran\c cois Berteloot, Universit\' e Paul Sabatier MIG. Laboratoire Emile Picard UMR 5580. F-31062 

Toulouse
Cedex 9 France.

{\em Email: berteloo@picard.ups-tlse.fr}}

\end{document}